\documentclass[11pt]{article}
\usepackage{mst-stylefile}
\usepackage{Bbb11a}
\usepackage{harvard}
\usepackage{amsfonts}
\usepackage{amssymb}
\usepackage{verbatim}
\include{epsf}

\catcode`@=12
\newcommand{\bbZ}{{\Bbb Z}}
\newcommand{\bbR}{{\Bbb R}}

\newcommand{\bbC}{{\Bbb C}}

\renewcommand{\cite}{\citeyear}

\begin{document}

\title{Identification of periodic and cyclic \\ fractional stable motions
\thanks{ This research was partially supported by the NSF grant DMS-0102410 at Boston University.}
\thanks{{\em AMS Subject classification}. Primary 60G18, 60G52;
secondary 28D, 37A.}
\thanks{{\em Keywords and phrases}: stable, self-similar
processes with stationary increments, mixed moving averages, dissipative and conservative
flows, periodic and cyclic flows, periodic and cyclic fractional stable motions.} }

\author{
 Vladas Pipiras
\\  University of North Carolina at Chapel Hill
\and
 Murad S.Taqqu
\\  Boston University
}

\bibliographystyle{agsm}

\maketitle



\begin{abstract}
\noindent Self-similar stable mixed moving average processes can be related to
nonsingular flows through their minimal representations. Self-similar stable mixed moving
averages related to dissipative flows have been studied, as well as processes associated
with identity flows which are the simplest type of conservative flows. The focus here is
on self-similar stable mixed moving averages related to periodic and cyclic flows.
Periodic flows are conservative flows such that each point in the space comes back to its
initial position in finite time, either positive or null. The flow is cyclic if the
return time is positive.

Self-similar mixed moving averages are called periodic, resp.\ cyclic, fractional stable
motions if their minimal representations are generated by periodic, resp.\ cyclic, flows.
These processes, however, are often defined by a nonminimal representation. We provide a
way to detect whether they are periodic or cyclic even if their representation is
nonminimal. By using these identification results, we obtain a more refined decomposition
of self-similar mixed moving averages.
\end{abstract}


\section{Introduction}
\label{s:preliminar}

Consider continuous-time stochastic processes $\{X(t)\}_{t\in\bbR}$ which have {\it
stationary increments} and are {\it self-similar} with self-similarity parameter $H>0$.
Stationarity of the increments means that the processes $X(t+h) - X(h)$ and $X(t) - X(0)$
have the same finite-dimensional distributions for any fixed $h\in\bbR$. Self-similarity
means that, for any fixed $c>0$, the processes $X(ct)$ and $c^H X(t)$ have the same
finite-dimensional distributions. The parameter $H>0$ is called the self-similarity
parameter. Self-similar stationary increments processes are of interest because their
increments can be used as models for stationary, possibly strongly dependent time series.

Fractional Brownian motion is the only (up to a multiplicative constant and for fixed
$H\in(0,1)$) Gaussian $H$--self-similar process with stationary increments. See, for
example, Section 7 in Samorodnitsky and Taqqu \cite{samorodnitsky:taqqu:1994book},
Embrechts and Maejima \cite{embrechts:maejima:2002} or two recent collections Doukhan,
Oppenheim and Taqqu \cite{doukhan:oppenheim:taqqu:2003} and Rangarajan and Ding
\cite{rangarajan:ding:2003} of survey articles. In contrast, for $\alpha\in(0,2)$, there
are infinitely many non-Gaussian $\alpha$-stable self-similar processes with stationary
increments. In Pipiras and Taqqu \cite{pipiras:taqqu:2002d,pipiras:taqqu:2002s}, the
authors have started to classify an important subclass of such processes, called {\it
self-similar mixed moving averages}, by relating them to ``flows'', an idea which has
originated with Rosi{\' n}ski \cite{rosinski:1995}. In this paper, we focus on
self-similar mixed moving averages which are related to {\it periodic} and, more
specifically, {\it cyclic} flows in the sense of Pipiras and Taqqu
\cite{pipiras:taqqu:2002d,pipiras:taqqu:2002s}. We call such processes {\it periodic} and
{\it cyclic fractional stable motions}. We show how, given a representation of the
process, one can determine whether a general self-similar mixed moving average is a
periodic or cyclic fractional stable motion. This leads to a decomposition of
self-similar mixed moving averages which is more refined than that obtained in Pipiras
and Taqqu \cite{pipiras:taqqu:2002d,pipiras:taqqu:2002s}. In a subsequent paper Pipiras
and Taqqu \cite{pipiras:taqqu:2003ca}, we study the properties of periodic and cyclic
fractional stable motions in greater detail, provide examples and show that periodic
fractional stable motions have canonical representations.

Many ideas are adapted from Pipiras and Taqqu \cite{pipiras:taqqu:2003cy} where we
investigated stable {\it stationary} processes related to periodic and cyclic flows in
the sense of Rosi{\' n}ski \cite{rosinski:1995}. Since these ideas appear in a simpler
form in Pipiras and Taqqu \cite{pipiras:taqqu:2003cy}, we suggest that the reader refers
to that paper for further clarifications and insight. The focus here is on stationary
increments mixed moving averages which are self-similar. Their connection to flows is
more involved and the results obtained in the stationary case cannot be readily applied.

Our presentation is different from that of Pipiras and Taqqu \cite{pipiras:taqqu:2003cy}.
While in Pipiras and Taqqu \cite{pipiras:taqqu:2003cy}, we focused first on stationary
stable processes having an {\it arbitrary} representation, we focus first here on
periodic and cyclic fractional stable motions having a ``minimal representation''. It is
convenient to work first with minimal representations because periodic and cyclic
fractional motions with minimal representations can be directly related to periodic and
cyclic flows. We then turn to self-similar mixed moving averages having an arbitrary,
possibly nonminimal, representation. This approach sheds additional light on the various
relations between stable processes and flows, and their corresponding decompositions in
disjoint classes.

Self-similar stable mixed moving averages, flows, cocycles and semi-additive functionals
are introduced in the next section. In Section \ref{s:intro}, we describe our results and
outline the rest of the paper.


\section{Description of the results}
\label{s:intro}

We shall focus on symmetric $\alpha$-stable ($S\alpha S$, in short), $\alpha \in (0,2)$,
self-similar processes $\{X_\alpha(t)\}_{t\in\bbR}$ with a {\it mixed moving average}
representation
\begin{equation}\label{e:mma}
\{X_\alpha(t)\}_{t\in\bbR} \stackrel{d}{=}\left\{ \int_X\int_\bbR \Big(G(x,t+u) -
G(x,u)\Big) M_\alpha(dx,du)\right\}_{t\in\bbR},
\end{equation}
where $\stackrel{d}{=}$ stands for the equality in the sense of the finite-dimensional
distributions. Here, $(X,{\cal X},\mu)$ is a standard Lebesgue space, that is, $(X,{\cal
X})$ is a measurable space with one-to-one, onto and bimeasurable correspondence to a
Borel subset of a complete separable metric space, and $\mu$ is a $\sigma$-finite
measure. $M_\alpha$ is a $S\alpha S$ random measure on $X\times \bbR$ with the control
measure $\mu(dx)du$ and $G:X\times\bbR\mapsto \bbR$ is some measurable deterministic
function. Saying that the process $X_\alpha$ is given by the representation (\ref{e:mma})
is equivalent to having its characteristic function expressed as
\begin{equation}\label{e:char-func-mma}
E\exp\Big\{i\sum_{k=1}^n \theta_k X_\alpha(t_k) \Big\} = \exp\biggl\{-\int_X \int_\bbR
\Big| \sum_{k=1}^n \theta_k G_{t_k}(x,u) \Big|^\alpha \mu(dx)du \biggr\},
\end{equation}
where
\begin{equation}\label{e:G_t}
G_t(x,u) = G(x,t+u) - G(x,u), \quad x\in X,u\in\bbR,
\end{equation}
and $\{G_t\}_{t\in\bbR} \subset L^\alpha(X\times\bbR,\mu(dx)du)$. The function
$G_t(x,u)$, or sometimes the function $G$, is called a {\it kernel function} of the
representation (\ref{e:mma}). For more information on $S\alpha S$ random measures and
integral representations of the type \refeq{mma}, see for example Samorodnitsky and Taqqu
\cite{samorodnitsky:taqqu:1994book}. Moreover, by setting $\xi = \sum_{k=1}^n \theta_k
X_\alpha(t_k)$, relation \refeq{char-func-mma} implies that $E \exp\{i\theta \xi\} =
\exp\{-\sigma^\alpha |\theta|^\alpha\}$ for some $\sigma\geq 0$ and all $\theta\in\bbR$.
By definition, $\xi$ is a $S\alpha S$ random variable and hence $X_\alpha$ is a $S\alpha
S$ process as well.

It follows from (\ref{e:char-func-mma}) that a mixed moving average $X_\alpha$ has always
stationary increments. Additional assumptions have to be imposed on the function $G$ for
the process $X_\alpha$ to be also self-similar. These assumptions are stated in
Definition \ref{d:ss-mma} and are formulated in terms of flows and some additional
functionals which we now define (see also Pipiras and Taqqu \cite{pipiras:taqqu:2002d}).

A (multiplicative) {\it flow} $\{\psi_c\}_{c>0}$ on $(X,{\cal X},\mu)$ is a collection of
deterministic measurable maps $\psi_c: X \to X$ satisfying
\begin{equation}\label{e:flow0}
\psi_{c_1c_2}(x) = \psi_{c_1}(\psi_{c_2}(x)), \quad \mbox{for all} \ c_1,c_2>0,\ x\in X,
\end{equation}
and $\psi_1(x)=x$ for all $x \in X$. The flow is {\it nonsingular} if each map $\psi_c,
c>0$, is nonsingular, that is, $\mu(A) = 0$ implies $\mu(\psi_c^{-1}(A))=0$. It is {\it
measurable} if a map $\psi_c(x): (0, \infty) \times X \to X$ is measurable.

A {\it cocycle} $\{b_c\}_{c>0}$ for the flow $\{\psi_c\}_{c>0}$ taking values in
$\{-1,1\}$ is a measurable map $b_c(x):(0, \infty) \times X \to \{-1,1\}$ satisfying
\begin{equation}\label{e:cocycle0}
b_{c_1c_2}(x) = b_{c_1}(x) b_{c_2}(\psi_{c_1}(x)), \quad \mbox{for all} \ c_1,c_2>0,\
x\in X.
\end{equation}

A {\it semi-additive functional} $\{g_c\}_{c>0}$ for the flow $\{\psi_c\}_{c>0}$ is a
measurable map $g_c(x):(0, \infty) \times X \to \bbR$ such that
\begin{equation}\label{e:1-semi-add0}
g_{c_1c_2}(x) = c_2^{-1} g_{c_1}(x) + g_{c_2}(\psi_{c_1}(x)), \quad \mbox{for all} \
c_1,c_2>0,\ x\in X.
\end{equation}

We use throughout the paper the useful notation
\begin{equation}\label{e:kappa}
    \kappa = H - \frac{1}{\alpha}.
\end{equation}
The support of $\{f_t\}_{t\in \bbR}\subset L^0(S,{\cal S},m)$, denoted
$\mbox{supp}\{f_t,t\in \bbR\}$, is a minimal (a.e.) set $A\in {\cal S}$ such that
$m\{f_t(s)\neq 0,s\notin A\}=0$ for every $t\in \bbR$.

\begin{definition}\label{d:ss-mma}
A $S\alpha S$, $\alpha \in (0,2)$, self-similar process $X_\alpha$ having a mixed moving
average representation (\ref{e:mma}) is said to be {\it generated by a nonsingular
measurable flow} $\{\psi_c\}_{c>0}$ on $(X,{\cal X},\mu)$ (through the kernel function
$G$) if

$(i)$ for all $c>0$,
\begin{equation}\label{e:generated}
c^{-\kappa} G(x,cu) = b_c(x) \left\{{d(\mu\circ \psi_c) \over d\mu}
(x)\right\}^{1/\alpha} G\Big(\psi_c(x), u + g_c(x)\Big) + j_c(x)\quad \mbox{a.e.}\
\mu(dx)du,
\end{equation}
where $\{b_c\}_{c>0}$ is a cocycle (for the flow $\{\psi_c\}_{c>0}$) taking values in
$\{-1,1\}$, $\{g_c\}_{c>0}$ is a semi-additive functional (for the flow
$\{\psi_c\}_{c>0}$) and $j_c(x)$ is some function, and

$(ii)$
\begin{equation}\label{e:full-support}
\mbox{supp}\left\{ G(x,t+u) - G(x,u),t\in\bbR\right\} = X\times\bbR \quad \mbox{a.e.}\
\mu(dx)du.
\end{equation}
\end{definition}

\bigskip
Relation (\ref{e:full-support}) is imposed in order to eliminate ambiguities stemming
from taking too big a space $X$. Definition \ref{d:ss-mma} can be found in Pipiras and
Taqqu \cite{pipiras:taqqu:2002d}. Observe that it involves the kernel $G$ and hence the
representation (\ref{e:mma}) of $X_\alpha$. The process $X_\alpha$ may have equivalent
representations (in the sense of the finite-dimensional distributions), each involving a
different function $G$. The so-called ``minimal representations'' are of particular
interest. Minimal representations were introduced by Hardin \cite{hardin:1982b} and
subsequently developed by Rosi\'nski \cite{rosinski:1998}. See also Section 4 in Pipiras
and Taqqu \cite{pipiras:taqqu:2002d}, or Appendix B in Pipiras and Taqqu
\cite{pipiras:taqqu:2003cy}. The representation $\{G_t\}_{t \in \bbR}$ of (\ref{e:mma})
is {\it minimal} if (\ref{e:full-support}) holds, and if for any nonsingular map $\Phi: X
\times \bbR \to X \times \bbR$ such that, for any $t \in \bbR$,
\begin{equation} \label{e:G_t-minimal}
G_t(\Phi(x,u))=k(x,u) G_t(x,u) \quad \mbox{a.e.} \ \mu(dx)du
\end{equation}
with some $k(x,u) \neq 0$, we have $\Phi(x,u)=(x,u)$, that is, $\Phi$ is the identity
map, a.e.\ $\mu(dx)du$.

Definition \ref{d:ss-mma} is closely related to self-similarity. By using
(\ref{e:char-func-mma}) together with (\ref{e:flow0})--(\ref{e:1-semi-add0}), it is easy
to verify that a mixed moving average (\ref{e:mma}) with a function $G$ satisfying
(\ref{e:generated}) is self-similar (see Pipiras and Taqqu \cite{pipiras:taqqu:2002d}).
Conversely, by Theorems 4.1 and 4.2 of Pipiras and Taqqu \cite{pipiras:taqqu:2002d}, any
$S\alpha S$, $\alpha \in (1,2)$, a self-similar mixed moving average is generated by a
flow in the sense of Definition \ref{d:ss-mma} with the kernel $G$ in (\ref{e:generated})
associated with the minimal representation of the process.

By using the connection between processes and flows, we proved in Pipiras and Taqqu
\cite{pipiras:taqqu:2002d} that $S\alpha S$, $\alpha \in (1,2)$, self-similar mixed
moving averages can be decomposed uniquely (in distribution) into two independent
processes as
\begin{equation}\label{e:into-two-parts}
    X_\alpha \stackrel{d}{=} X_\alpha^{D} + X_\alpha^{C}.
\end{equation}
Here, $X_\alpha^{D}$ is a self-similar mixed moving average generated by a {\it
dissipative flow}. Informally, the flow $\{ \psi_c\}_{c>0}$ is dissipative when the
points $x$ and $\psi_c(x)$ move further apart as $c$ approaches $\infty$ ($\ln
c\to\infty$) or $c$ approaches $0$ ($\ln c\to -\infty$). An example of a dissipative flow
is $\psi_c(x) = x +\ln c$, $x\in\bbR$. Self-similar mixed moving average processes
generated by dissipative flows have a canonical representation (see Theorem 4.1 in
Pipiras and Taqqu \cite{pipiras:taqqu:2002s}) and are studied in detail in Pipiras and
Taqqu \cite{pipiras:taqqu:2003di}, where they are called {\it dilated fractional stable
motions}.

The process $X_\alpha^C$ in (\ref{e:into-two-parts}) is a self-similar mixed moving
average generated by a {\it conservative flow}. Conservative flows $\{\psi_c\}_{c>0}$ are
such that the points $x$ and $\psi_c(x)$ become arbitrarily close at infinitely many
values of $c$. An example of a conservative flow is $\psi_c(x) = xe^{i\ln c}$, $|x|=1,\,
x\in\bbC$ since $\psi_c(x) = x$ every time that $\ln c$ is a multiple of $2\pi$. Although
this example is elementary, the general structure of conservative flows is complex and,
in particular, more intricate than that of dissipative flows. Consequently, contrary to
the processes generated by dissipative flows, there is no simple canonical representation
of the self-similar mixed moving averages generated by conservative flows.

It is nevertheless possible to obtain a further decomposition of self-similar mixed
moving averages generated by conservative flows. As shown in Pipiras and Taqqu
\cite{pipiras:taqqu:2002s},
\begin{equation}\label{e:into-three-parts}
X_\alpha^C\stackrel{d}{=}X_\alpha^{F} +  X_\alpha^{C\setminus F},
\end{equation}
where the decomposition is unique in distribution and has independent components. The
processes $X_\alpha^{F}$ in the decomposition (\ref{e:into-three-parts}) are those
self-similar mixed moving averages that have a canonical representation (\ref{e:mma})
with the kernel function
\begin{equation} \label{e:mlfsm-kernel}
G(x,u) = \left\{
    \begin{array}{ll}
    F_1(x) u_+^\kappa + F_2(x) u_-^\kappa, & \kappa \ne 0, \\
    F_1(x) \ln|u| + F_2(x) 1_{(0,\infty)}(u), & \kappa = 0,
    \end{array}
\right.
\end{equation}
where $u_+ = \max\{0,u\}$, $u_- = \max\{0,-u\}$ and $F_1,F_2:Z\mapsto \bbR$ are some
functions. Thus,
\begin{equation}\label{e:mlfsm}
X_\alpha^F(t) \stackrel{d}{=} \left\{
    \begin{array}{ll}
    \int_X\int_\bbR \Big(F_1(x)((t+u)_+^\kappa - u_+^\kappa) + F_2(x)((t+u)_-^\kappa - u_-^\kappa)\Big) M_\alpha(dx,du),& \kappa \ne 0, \\
    \int_X\int_\bbR \Big(F_1(x)\ln\frac{|t+u|}{|u|} + F_2(x)1_{(-t,0)}(u)\Big) M_\alpha(dx,du),& \kappa =
    0.
    \end{array}\right.
\end{equation}
The processes (\ref{e:mlfsm}) are called {\it mixed linear fractional stable motions}
(mixed LFSM, in short) and are essentially generated by {\it identity flows}
(`essentially' will become clear in the sequel). An identity flow is the simplest type of
conservative flow, defined by $\psi_c(x) = x$ for all $c>0$, and the upperscript $F$ in
$X_\alpha^F$ refers to the fact that the points $x$ are {\it fixed} points under the
flow. The processes $X_\alpha^{C\setminus F}$ in (\ref{e:into-three-parts}) are
self-similar mixed moving averages generated by conservative flows but without the {\it
mixed LFSM component} (\ref{e:mlfsm}), that is, they cannot be represented in
distribution by a sum of two independent processes, one of which is a nondegenerate mixed
LFSM (\ref{e:mlfsm}).

Our goal is to obtain a more detailed decomposition of self-similar mixed moving
averages. We will show that there are independent self-similar mixed moving averages
$X_\alpha^L$ and $X_\alpha^{C\setminus P}$ such that
\begin{equation} \label{e:into-four-parts1}
  X_\alpha^{C \setminus F} \stackrel{d}{=} X_\alpha^{L} + X_\alpha^{C \setminus P}
\end{equation}
and hence, in view of (\ref{e:into-three-parts}),
\begin{equation} \label{e:into-four-parts2}
  X_\alpha^{C} \stackrel{d}{=} X_\alpha^{F} + X_\alpha^{L} + X_\alpha^{C \setminus P}
  =: X_\alpha^{P} + X_\alpha^{C \setminus P},
\end{equation}
where the decompositions (\ref{e:into-four-parts1}) and (\ref{e:into-four-parts2}) are
unique in distribution and have independent components. While the processes
$X_\alpha^{F}$ are essentially generated by identity flows, the process $X_\alpha^{P} =
X_\alpha^F + X_\alpha^L$ and the process $X_\alpha^{L}$ are essentially generated by {\it
periodic} and {\it cyclic} flows, respectively\footnote{The letters $D$ and $C$ are
associated with {\it D}issipative and {\it C}onservative flows, respectively. The letter
$F$ (``{\it F}ixed'') is associated with identity flows, the letter $P$ with {\it
P}eriodic flows and the letter $L$ with cyc{\it L}ic flows.}. Periodic flows are
conservative flows such that any points in the space comes back to its initial position
in a finite period of time. Identity flows are periodic flows with period zero. Cyclic
flows are periodic flows with positive period. Cyclic flows are probably the simplest
type of conservative flows after the identity flows.

These flows are defined as follows. Let $\{\psi_c\}_{c>0}$ be a measurable flow on a
standard Lebesgue space $(X,{\cal X},\mu)$. Consider the following subsets of $X$ induced
by the flow $\{\psi_c\}_{c>0}$:
\begin{eqnarray}
P & := & \{x:\exists\ p = p(x) \neq 1: \psi_p(x) = x\},\label{e:P} \\
F & := & \{x: \psi_c(x) = x\ \mbox{for all}\ c>0\}, \label{e:F} \\
L &:= & P \setminus F. \label{e:L}
\end{eqnarray}

\begin{definition}\label{d:periodic-fixed-cyclic-points}
The elements of $P$, $F$, $L$ are called the {\it periodic}, {\it fixed}  and {\it
cyclic} points of the flow $\{\psi_c\}_{c>0}$, respectively.
\end{definition}

By Lemma 2.1 in Pipiras and Taqqu \cite{pipiras:taqqu:2003cy}, the sets $P,L$ are
$\mu$-measurable (measurable with respect to the measure $\mu$) and the set $F$ is
(Borel) measurable.

\begin{definition}\label{d:periodic-cyclic1}
A measurable flow $\{\psi_c\}_{c>0}$ on $(X,{\cal X},\mu)$ is {\it periodic} if $X=P$
$\mu$-a.e., is {\it identity} if $X=F$ $\mu$-a.e., and it is {\it cyclic} if $X=L$
$\mu$-a.e.
\end{definition}

An alternative equivalent definition of a cyclic flow is given in Definition 2.2 of
Pipiras and Taqqu \cite{pipiras:taqqu:2003cy}, but it will not be used here.

The processes $X_\alpha^{P}$ and $X_\alpha^{L}$ in (\ref{e:into-four-parts2}) will be
called, respectively, {\it periodic fractional stable motions} and {\it cyclic fractional
stable motions}. We indicated above that the processes $X_\alpha^{F}$, $X_\alpha^{P}$ and
$X_\alpha^{L}$ are essentially determined by identity, periodic and cyclic flows,
respectively. By `essentially determined', we mean that if the processes $X_\alpha^{P}$
and $X_\alpha^{L}$ are given by their minimal representations, then they are necessarily
generated by periodic and cyclic flows, respectively, in the sense of Definition
\ref{d:ss-mma}. This terminology is not restrictive in the case $\alpha\in (1,2)$ because
mixed moving averages always have minimal representations (\ref{e:mma}) by Theorem 4.2 in
Pipiras and Taqqu \cite{pipiras:taqqu:2002d} and, according to Theorem 4.1 of that paper,
self-similar mixed moving averages given by a minimal representation (\ref{e:mma}) are
always generated by a unique flow in the sense of Definition
\ref{d:ss-mma}.\footnote{When $\alpha\in (0,1]$, we were able to prove Theorem 4.2
concerning existence of minimal representations for mixed moving averages only under
additional assumptions on the process (see Remark 6 following Theorem 4.2). To keep the
presentation simple, we do not introduce here these additional assumptions and hence
suppose in this paper that $\alpha\in (1,2)$, unless stated explicitly otherwise. } More
generally, when a self-similar mixed moving average $X_\alpha$ given by a minimal
representation (\ref{e:mma}), is generated by the flow, the processes $X_\alpha^{P}$ and
$X_\alpha^{L}$ in the decomposition (\ref{e:into-four-parts2}) can be defined by
replacing respectively the space $X$ in the integral representation (\ref{e:mma}) by $P$
and $L$, that is, the periodic and cyclic point sets of the generating flow.

Why are we referring to minimal representations? If one makes no restrictions on the form
of a representation (\ref{e:mma}), periodic and cyclic fractional stable motions can be
generated by flows other than periodic and cyclic, and the components $X_\alpha^{P}$ and
$X_\alpha^{L}$ in the decomposition (\ref{e:into-four-parts2}) may not be related to the
periodic and cyclic point sets of the underlying flow. An analogous phenomenon is also
associated with the decomposition (\ref{e:into-three-parts}). Since we would like to work
with an arbitrary (not necessarily minimal) representation (\ref{e:mma}), it is desirable
to be able to recognize periodic and cyclic fractional stable motions without relying on
minimal representations and flows. We shall therefore provide identification criteria
based on the (possibly nonminimal) kernel function $G$ in the representation
(\ref{e:mma}). These criteria allow one to obtain the decompositions
(\ref{e:into-four-parts1}) and (\ref{e:into-four-parts2}) when starting with an arbitrary
(possibly nonminimal) representation (\ref{e:mma}).

The paper is organized as follows. In Section \ref{s:pcfsm-minimal}, we establish the
decompositions (\ref{e:into-four-parts1}) and (\ref{e:into-four-parts2}) using
representations (\ref{e:mma}) that are minimal, and introduce periodic and cyclic
fractional stable motions. Criteria to identify periodic and cyclic fractional stable
motions through (possibly nonminimal) kernel functions $G$ are provided in Sections
\ref{s:pfsm-nonminimal} and \ref{s:cfsm-nonminimal}. The decompositions
(\ref{e:into-four-parts1}) and (\ref{e:into-four-parts2}) which are based on these
criteria can be found in Section \ref{s:refined-decomposition}.


\section{Periodic and cyclic fractional stable motions: the minimal case}
\label{s:pcfsm-minimal}

By Theorem 4.2 in Pipiras and Taqqu \cite{pipiras:taqqu:2002d}, any $S \alpha S$,
$\alpha\in (1,2)$, mixed moving average $X_\alpha$ has an integral representation
(\ref{e:mma}) which is minimal. By Theorem 4.1 in Pipiras and Taqqu
\cite{pipiras:taqqu:2002d}, a self-similar mixed moving averages $X_\alpha$ given by a
minimal representation (\ref{e:mma}) is generated by a unique flow $\{\psi_c\}_{c>0}$ in
the sense of Definition \ref{d:ss-mma}.

By the Hopf decomposition (see e.g.\ Pipiras and Taqqu \cite{pipiras:taqqu:2002d}), the
space $X$ can be decomposed into two parts, $D$ and $C$, invariant under the flow, $D$
denoting the dissipative points of $\{\psi_c\}_{c>0}$ and $C$ denoting the conservative
points of $\{\psi_c\}_{c>0}$. Let $D$, $C$, $F$, $L$ and $P$ be then the dissipative,
conservative, fixed, cyclic and periodic point sets of the flow $\{\psi_c\}_{c>0}$,
respectively. Since
$$
X=D+C=D+P+C\setminus P = D+F+L+C\setminus P,
$$
we can write
\begin{equation}\label{e:into-four-minimal}
X_\alpha \stackrel{d}{=}  X_\alpha^D + X_\alpha^P + X_\alpha^{C\setminus P} = X_\alpha^D
+ X_\alpha^F + X_\alpha^L + X_\alpha^{C\setminus P},
\end{equation}
where
$$
X_\alpha^P = X_\alpha^F + X_\alpha^L
$$
and where, for a set $S\subset X$,
\begin{equation}\label{e:X_alpha^S}
X_\alpha^S(t) = \int_S\int_\bbR G_t(x,u) M_\alpha(dx,du).
\end{equation}
Since  by their definitions, the sets $D$, $C$, $F$, $P$ and $L$ are invariant under the
flow, the processes $X_\alpha^D$, $X_\alpha^F$, $X_\alpha^L$ and $X_\alpha^{C\setminus
P}$ are self-similar mixed moving averages. These processes are independent because the
sets $D$, $F$, $L$ and $C\setminus P$ are disjoint (see Theorem 3.5.3 in Samorodnitsky
and Taqqu \cite{samorodnitsky:taqqu:1994book}). The processes $X_\alpha^S$ are generated
by the flow $\psi^S$ where $\psi^S$ denotes the flow $\psi$ restricted to a set $S$,
which is invariant under the flow. Observe that $\psi^D$, $\psi^F$, $\psi^L$ and $\psi^P$
are dissipative, identity, cyclic and periodic flows, respectively, and that
$\psi^{C\setminus P}$ is a conservative flow without periodic points, and, for example,
the process $X_\alpha^D$ is generated by the dissipative flow $\psi^D$.

A self-similar mixed moving average may have another minimal representation (\ref{e:mma})
with a kernel function $\widetilde G$ on the space $\widetilde X$, and hence be generated
by another flow $\{\widetilde \psi_c\}_{c>0}$. Partitioning $\widetilde X$ into the
dissipative, fixed, cyclic and ``other'' conservative point sets of the flow
$\{\widetilde \psi_c\}_{c>0}$ as above, leads to the decomposition
\begin{equation}\label{e:into-four-minimal2}
X_\alpha \stackrel{d}{=}  \widetilde X_\alpha^{D} + \widetilde X_\alpha^{F} + \widetilde
X_\alpha^{L} + \widetilde X_\alpha^{C\setminus P}.
\end{equation}
We will say that the decomposition (\ref{e:into-four-minimal}) obtained from a minimal
representation (\ref{e:mma}) is {\it unique in distribution} if the distribution of its
components does not depend on the minimal representation used in the decomposition. In
other words, uniqueness in distribution holds if
\begin{equation}\label{e:unique-minimal}
X_\alpha^D \stackrel{d}{=} \widetilde X_\alpha^D,\quad X_\alpha^F \stackrel{d}{=}
\widetilde X_\alpha^F,\quad X_\alpha^L \stackrel{d}{=} \widetilde X_\alpha^L,\quad
X_\alpha^{C\setminus P} \stackrel{d}{=} \widetilde X_\alpha^{C\setminus P},
\end{equation}
where $X_\alpha^S$ and $\widetilde X_\alpha^S$ with $S=D,F,L$ and $C\setminus P$, are the
components of the decompositions (\ref{e:into-four-minimal}) and
(\ref{e:into-four-minimal2}) obtained from two different minimal representations of the
process.

\begin{theorem}\label{t:unique-minimal}
Let $\alpha\in (1,2)$. The decomposition (\ref{e:into-four-minimal}) obtained from a
minimal representation (\ref{e:mma}) of a self-similar mixed moving average $X_\alpha$ is
unique in distribution.
\end{theorem}

\begin{proof}
Suppose that a self-similar mixed moving average $X_\alpha$ is given by two different
minimal representations with the kernel functions $G$ and $\widetilde G$, and the spaces
$(X,\mu)$ and $(\widetilde X,\widetilde \mu)$, respectively. Suppose also that $X_\alpha$
is generated through these minimal representations by two different flows
$\{\psi_c\}_{c>0}$ and $\{\widetilde \psi_c\}_{c>0}$ on the spaces $X$ and $\widetilde
X$, respectively. Let (\ref{e:into-four-minimal}) and (\ref{e:into-four-minimal2}) be two
decompositions of $X_\alpha$ obtained from these two minimal representations and the
generating flows. Let also $D,F,L,P,C$ and $\widetilde D,\widetilde F,\widetilde
L,\widetilde P,\widetilde C$ be the dissipative, fixed, cyclic, periodic and conservative
point sets of the flows $\{\psi_c\}_{c>0}$ and $\{\widetilde \psi_c\}_{c>0}$,
respectively. We need to show that the equalities (\ref{e:unique-minimal}) hold.

By Theorem 4.3 and its proof in Pipiras and Taqqu \cite{pipiras:taqqu:2002d}, the kernel
functions $G$ and $\widetilde G$, and the flows $\psi$ and $\widetilde\psi$ are related
in the following way. There is a map $\Phi:\widetilde X\mapsto X$ such that $(i)$ $\Phi$
is one-to-one, onto and bimeasurable (up to two sets of measure zero); $(ii)$ $\widetilde
\mu\circ \Phi$ and $\mu$ are mutually absolutely continuous; $(iii)$ for all $c>0$,
$\psi_c\circ \Phi = \Phi \circ \widetilde \psi_c$\, $\widetilde\mu$-a.e., and $(iv)$ for
all $t\in\bbR$,
\begin{equation}\label{e:G-tildeG-minimal}
\widetilde G_t(\widetilde x,u)=b(\widetilde x)\left\{\frac{d(\mu\circ \Phi)}{d\widetilde
\mu}(\widetilde x)\right\}^{1/\alpha} G_t(\Phi(\widetilde x),u+g(\widetilde x)),\quad
\mbox{a.e.}\ \widetilde \mu(d\widetilde x)du,
\end{equation}
where $b:\widetilde X\mapsto\{-1,1\}$ and $g:\widetilde X\mapsto \bbR$ are measurable
functions.

Since $D$ ($C$, resp.) can be expressed as
$$
D\ (C,\mbox{resp.}) = \left\{x\in X:\int_0^\infty f(\psi_c(x)) \frac{d(\mu\circ
\psi_c)}{d\mu}(x)\, c^{-1}dc <\infty \ (=\infty,\ \mbox{resp.}) \right\},\quad
\mu-\mbox{a.e.},
$$
for any $f\in L^1(X,\mu)$, $f>0$ a.e.\ (see, for example, (3.22) and (3.33) in Pipiras
and Taqqu \cite{pipiras:taqqu:2002d} in the case of additive flows), we obtain by using
the relations $(ii)$ and $(iii)$ above that
\begin{equation}\label{e:DC-tildeDC-minimal}
\Phi^{-1}(D) = \widetilde D,\quad \Phi^{-1}(C) = \widetilde C, \quad \widetilde
\mu-\mbox{a.e.}
\end{equation}
By using relations $(i)$--$(iii)$, we can deduce directly from (\ref{e:F})--(\ref{e:L})
that
\begin{equation}\label{e:FL-tildeFL-minimal}
\Phi^{-1}(F) = \widetilde F,\quad \Phi^{-1}(P) = \widetilde P,\quad \Phi^{-1}(L) =
\widetilde L, \quad \widetilde \mu\mbox{-a.e.}
\end{equation}
and hence
\begin{equation}\label{e:PC-tildePC-minimal}
\Phi^{-1}(C\setminus P) = \widetilde C\setminus \widetilde P, \quad \widetilde
\mu\mbox{-a.e.}
\end{equation}

The equalities (\ref{e:unique-minimal}) can now be obtained by using
(\ref{e:G-tildeG-minimal}) together with
(\ref{e:DC-tildeDC-minimal})--(\ref{e:PC-tildePC-minimal}). For example, the first
equality in (\ref{e:unique-minimal}) follows by using (\ref{e:G-tildeG-minimal}) and
(\ref{e:DC-tildeDC-minimal}) to show that
$$
\int_{\widetilde D}\int_\bbR \left| \sum_{k=1}^n \theta_k\Big( \widetilde G(\widetilde
x,t_k+u) - \widetilde G(\widetilde x,u) \Big)\right|^\alpha \widetilde\mu(d\widetilde
x)du
$$
$$
=\int_{\Phi^{-1}(D)}\int_\bbR \left| \sum_{k=1}^n \theta_k\Big( G(\Phi(\widetilde
x),t_k+u+g(\widetilde x)) - G(\Phi(\widetilde x),u+g(\widetilde x)) \Big)\right|^\alpha
\frac{d(\mu\circ\Phi)}{d\widetilde \mu}(\widetilde x)\, \widetilde\mu(d\widetilde x)du
$$
$$
= \int_{\Phi^{-1}(D)}\int_\bbR \left| \sum_{k=1}^n \theta_k\Big( G(\Phi(\widetilde
x),t_k+u) - G(\Phi(\widetilde x),u) \Big)\right|^\alpha (\mu\circ\Phi)(d\widetilde x)du
$$
$$
= \int_{D}\int_\bbR \left| \sum_{k=1}^n \theta_k\Big( G(x,t_k+u) - G(x,u)
\Big)\right|^\alpha \mu(dx)du,
$$
where in the last equality, we used a change of variables. \ \ $\Box$
\end{proof}

\bigskip

Since the decomposition (\ref{e:into-four-minimal}) can be obtained through a minimal
representation for any $S\alpha S$, $\alpha\in(1,2)$, self-similar mixed moving average,
and it is unique in distribution by Theorem \ref{t:unique-minimal}, we may give the
following definition.

\begin{definition}\label{d:pfsm-cfsm-minimal}
A $S\alpha S$, $\alpha\in(1,2)$, self-similar mixed moving average $X_\alpha$ is called
{\it periodic fractional stable motion} ({\it cyclic fractional stable motion}, resp.) if
$$
X_\alpha \stackrel{d}{=} X_\alpha^P\quad (X_\alpha \stackrel{d}{=} X_\alpha^L,\
\mbox{resp.}),
$$
where $X_\alpha^L$ and $X_\alpha^P$ are the two components in the decomposition
(\ref{e:into-four-minimal}) of $X_\alpha$ obtained through a minimal representation.
\end{definition}

\noindent {\bf Notation.} Periodic and cyclic fractional stable motion will be
abbreviated as {\it PFSM} and {\it CFSM}, respectively.

\medskip
An equivalent definition of periodic and cyclic fractional stable motions is as follows.

\begin{proposition}\label{p:pfsm-cfsm-minimal-flow}
A $S\alpha S$, $\alpha\in (1,2)$, self-similar mixed moving average is a periodic
(cyclic, resp.) fractional stable motion if and only if the generating flow corresponding
to its minimal representation is periodic (cyclic, resp.).
\end{proposition}

\begin{proof}
By Definition \ref{d:pfsm-cfsm-minimal}, a self-similar mixed moving average $X_\alpha$
is a PFSM (CFSM, resp.) if and only if $X_\alpha =_d X_\alpha^P$ ($X_\alpha =_d
X_\alpha^L$, resp.), where $P$ ($L$, resp.) is the set of periodic (cyclic, resp.) points
of the generating flow $\psi$ corresponding to a minimal representation. It follows from
(\ref{e:into-four-minimal}) and (\ref{e:X_alpha^S}) that $X_\alpha =_d X_\alpha^P$
($X_\alpha =_d X_\alpha^L$, resp.) if and only if $X=P$ ($X=L$, resp.) $\mu$-a.e.\ and
hence, by Definition \ref{d:periodic-cyclic1}, if and only if the flow $\psi$ is periodic
(cyclic, resp.). \ \ $\Box$
\end{proof}

\bigskip

Definition \ref{d:pfsm-cfsm-minimal} and Proposition \ref{p:pfsm-cfsm-minimal-flow} use
minimal representations. Minimal representations, however, are not very easy to determine
in practice. It is therefore desirable to recognize a PFSM and a CFSM based on any,
possibly nonminimal representation. Since many self-similar mixed moving averages given
by nonminimal representations are generated by a flow in the sense of Definition
\ref{d:ss-mma}, we could expect that the process is a PFSM (CFSM, resp.) if the
generating flow is periodic (cyclic, resp.). This, however, is not the case in general.
For example, if a PFSM or CFSM $X_{\alpha}(t)= \int_X \int_{\bbR} G_t(x,u) M_{\alpha}
(dx,du)$ is generated by a periodic or cyclic flow $\psi_c(x)$ on $X$, we can also
represent the process $X_{\alpha}$ as $\int_Y \int_X \int_{\bbR}  G_t(x,u) M_{\alpha}
(dy,dx,du)$, where $ G_t(x,u)$ does not depend on $y$ and the control measure $\eta(dy)$
of $M_{\alpha} (dy,dx,du)$ in the variable $y$ is such that $\eta(Y)=1$. If $\widetilde
\psi_c(y)$ is a measure preserving flow on $(Y, \eta)$, then the process $X_{\alpha}$ is
also generated by the flow $\Phi_c(y,x) = (\widetilde \psi_c(y), \psi_c(x))$ on $Y \times
X$. If, in addition, the flow $\widetilde \psi_c(y)$ is not periodic (and hence not
cyclic), then the flow $\Phi_c(y,x)$ is neither periodic nor cyclic.

We will provide identification criteria for a PFSM and a CFSM which do not rely on either
minimal representations or flows, and which are based instead on the structure of the
kernel function $G$. An analogous approach was taken by Rosi{\' n}ski
\cite{rosinski:1995} to identify harmonizable processes, by Pipiras and Taqqu
\cite{pipiras:taqqu:2002s} to identify a mixed LFSM, and by Pipiras and Taqqu
\cite{pipiras:taqqu:2003cy} to identify periodic and cyclic stable stationary processes.


\section{Identification of periodic fractional stable motions: the nonminimal
case}\label{s:pfsm-nonminimal}

We first provide a criterion to identify periodic fractional stable motions without using
flows or minimal representations. The criterion is based on the periodic fractional
stable motion set which we define next. Let $X_\alpha$ be a self-similar mixed moving
average (\ref{e:mma}) defined through a (possibly nonminimal) kernel function $G$.

\begin{definition}\label{d:pfsm-set}
A {\it periodic fractional stable motion set} ({\it PFSM set}, in short) of a
self-similar mixed moving average $X_\alpha$ given by (\ref{e:mma}), is defined as
\begin{eqnarray}
C_P & = & \Big\{ x\in X: \exists\  c = c(x)\neq 1: G(x,cu) = b\, G(x,u + a)
+ d\ \ \mbox{a.e.}\ du\nonumber \\
& &\hspace{.7in} \mbox{for some}\ a = a(c,x), b=b(c,x)\neq 0,d=d(c,x)\in\bbR \Big\}.
\label{e:pfsm-set}
\end{eqnarray}
\end{definition}

\begin{proposition}\label{p:pfsm-set-alt}
The relation in (\ref{e:pfsm-set}) can be expressed as
\begin{equation}\label{e:pfsm-set-exp1}
G(x,cu+g) = b\, G(x,u + g) + d
\end{equation}
for some $b\neq 0$, $c\neq 1$, $g,d\in\bbR$. When $b\neq 1$, it can also be expressed as
\begin{equation}\label{e:pfsm-set-exp2}
G(x,cu+g) + f = b (G(x,u+g) + f)
\end{equation}
for some $b\neq 0$, $c\neq 1$, $g,f\in\bbR$.
\end{proposition}

\begin{proof}
Relation (\ref{e:pfsm-set-exp1}) follows by making the change of variables $u = v +
a/(c-1)$ in $G(x,cu) = b G(x,u + a) + d$. When $b\neq 1$, by writing $d=bf-f$ with
$f=d/(b-1)$ in (\ref{e:pfsm-set-exp1}), we get (\ref{e:pfsm-set-exp2}). \ \ $\Box$
\end{proof}

\medskip
Whereas the set of periodic points $P$ is defined by (\ref{e:P}) in terms of the flow
$\{\psi_c\}_{c>0}$, the set $C_P$ in (\ref{e:pfsm-set}) is defined in terms of the kernel
$G$. Definition \ref{d:pfsm-set} states that there is a factor $c$ such that the kernel
$G$ at time $u$ is related to the kernel at time $cu$.

\begin{lemma}\label{l:measurability-of-L}
The PFSM set $C_P$ in (\ref{e:pfsm-set}) is $\mu$-measurable. Moreover, the functions
$c(x),a(x) = a(c(x),x)$, $b=b(c(x),x)$ and $d=d(c(x),x)$ in (\ref{e:pfsm-set}) can be
taken to be $\mu$-measurable as well.
\end{lemma}

\begin{proof}
We first show the measurability of $C_L$. Consider the set
$$
A = \Big\{ (x,c,a,b,d) : G(x,cu) = b\, G(x,u+a) + d\ \ \mbox{a.e.} \ du\Big\}.
$$
Since $A = \{F(x,c,a,b,d) = 0\}$, where the function
$$
F(x,c,a,b,d) = \int_\bbR 1_{\{G(x,cu) = b\, G(x,u + a) + d \}}(x,c,a,b,d,u) du
$$
is measurable by the Fubini's theorem, we obtain that the set $A$ is measurable. Observe
that the set $C_P$ is a projection of the set $A$ on $x$, namely, that
$$
C_P = \mbox{proj}_X A : = \{x:\exists\ c,a,b,d: (x,c,a,b,d) \in A\}.
$$
Lemma 3.3 in Pipiras and Taqqu \cite{pipiras:taqqu:2003cy} implies that the PFSM set
$C_P$ is $\mu$-measurable and that the functions $a(x)$, $b(x)$, $c(x)$ and $d(x)$ can be
taken to be $\mu$-measurable as well.\ \ $\Box$
\end{proof}

\bigskip
In the next theorem, we characterize a PFSM in terms of the set $C_P$ instead of using
the set $P$ which involves flows as is done in Definition \ref{d:pfsm-set} and
Proposition \ref{p:pfsm-cfsm-minimal-flow}. Flows and minimal representations, however,
are used in the proof.

\begin{theorem}\label{t:identification-pfsm}
A $S \alpha S$, $\alpha\in (1,2)$, self-similar mixed moving average $X_{\alpha}$ given
by (\ref{e:mma}) with $G$ satisfying (\ref{e:full-support}) is a PFSM if and only if $C_P
= X$ $\mu$-{\rm a.e.}, where $C_P$ is the PFSM set defined in (\ref{e:pfsm-set}).
\end{theorem}

\begin{proof} Suppose first that $X_\alpha$ is a self-similar mixed moving average
given by (\ref{e:mma}) with $G$ satisfying (\ref{e:full-support}) and such that $C_P = X$
$\mu$-a.e. To show that $X_\alpha$ is a PFSM, we adapt the proof of Theorem 3.2 in
Pipiras and Taqqu \cite{pipiras:taqqu:2003cy}. The proof consists of 2 steps.

\medskip
{\it Step 1:} We will show that without loss of generality, the representation
(\ref{e:mma}) can be supposed to be minimal with $C_P = X$ $\mu$-a.e.\ By Theorem 4.2 in
Pipiras and Taqqu \cite{pipiras:taqqu:2002d}, the process $X_\alpha$ has a minimal
integral representation
\begin{equation}\label{e:minimal-mma}
    \int_{\widetilde X}\int_\bbR \Big( \widetilde G(\widetilde x,t+u) -
    \widetilde G(\widetilde x,u)\Big) \widetilde M_\alpha(d\widetilde x,du),
\end{equation}
where $(\widetilde X,\widetilde{\cal X},\widetilde \mu)$ is a standard Lebesgue space and
$\widetilde M_\alpha(d\widetilde x,du)$ has control measure $\widetilde \mu(d\widetilde
x) du$. Letting $\widetilde C_P$ be the periodic component set of $X_\alpha$ defined
using the kernel function $\widetilde G$, we need to show that $\widetilde C_P =
\widetilde X$ $\widetilde\mu$-a.e. By Corollary 5.1 in Pipiras and Taqqu
\cite{pipiras:taqqu:2002d}, there are measurable maps $\Phi_1:X\mapsto \widetilde X$,
$h:X\mapsto \bbR\setminus\{0\}$ and $\Phi_2,\Phi_3:X\mapsto \bbR$ such that
\begin{equation}\label{e:nonminimal-minimal}
    G(x,u) = h(x) \widetilde G(\Phi_1(x),u+\Phi_2(x)) + \Phi_3(x)
\end{equation}
a.e.\ $\mu(dx)du$, and
\begin{equation}\label{e:nonminimal-mu-minimal-mu}
    \widetilde \mu = \mu_h\circ \Phi_1^{-1},
\end{equation}
where $\mu_h(dx) = |h(x)|^\alpha \mu(dx)$. If $x\in C_P$, then
\begin{equation}\label{e:generated-c(x)}
G(x,c(x)u) = b(x)G(x,u+a(x)) + d(x)\quad \mbox{a.e.}\ du,
\end{equation}
for some functions $a(x),b(x)$, $c(x)$ and $d(x)$. Hence, by using
(\ref{e:nonminimal-minimal}) and (\ref{e:generated-c(x)}), we have for some functions
$F_1$,$F_2$ and $F_3$, a.e.\ $\mu(dx)$,
$$
\widetilde G(\Phi_1(x),c(x)u+\Phi_2(x)) = (h(x))^{-1} G(x,c(x)u) + F_1(x)
$$
$$
= (h(x))^{-1}b(x)G(x,u+a(x)) + F_2(x) = b(x) \widetilde G(\Phi_1(x),u+a(x)+\Phi_2(x)) +
F_3(x)
$$
a.e.\ $du$. This shows that $\Phi_1(x)\in \widetilde C_P$ and hence
\begin{equation}\label{e:C_P-subset-Phi^-1-tile-C_P}
C_P\subset \Phi_1^{-1}(\widetilde C_P),\quad \mu\mbox{-a.e.}
\end{equation}
Since $C_P = X$ $\mu$-a.e., we have $X = \Phi^{-1}(\widetilde C_P)$ $\mu$-a.e. This
implies $\widetilde C_P = \widetilde X$ $\widetilde \mu$-a.e., because if $\widetilde
\mu(\widetilde X\setminus \widetilde C_P)>0$, then by (\ref{e:nonminimal-mu-minimal-mu}),
we have $\mu(\Phi_1^{-1}(\widetilde X\setminus \widetilde C_P)) =
\mu(\Phi_1^{-1}(\widetilde X)\setminus X) = \mu(\emptyset) >0$.

\medskip
{\it Remark:} The converse is shown in the same way: if $C_P$ is not equal to $X$
$\mu$-a.e., then $\Phi_1^{-1}(\widetilde C_P)\subset C_P$ $\mu$-a.e. Together with
(\ref{e:C_P-subset-Phi^-1-tile-C_P}), this implies
\begin{equation}\label{e:C_P-is-Phi^-1-tile-C_P}
C_P= \Phi_1^{-1}(\widetilde C_P),\quad \mu\mbox{-a.e.}
\end{equation}
The relation (\ref{e:C_P-is-Phi^-1-tile-C_P}) is used in the proof of the converse of
this theorem and in the proof of Theorem \ref{t:into-four} below.

\medskip
We may therefore suppose without loss of generality that the representation (\ref{e:mma})
is minimal and that $C_P = X$ $\mu$-a.e. By Theorem 4.1 in Pipiras and Taqqu
\cite{pipiras:taqqu:2002d}, since the representation (\ref{e:mma}) is minimal, the
process $X_\alpha$ is generated by a flow $\{\psi_c\}_{c>0}$ and related functionals
$\{b_c\}_{c>0}$, $\{g_c\}_{c>0}$ and $\{j_c\}_{c>0}$ in the sense of Definition
\ref{d:ss-mma}.

\medskip
{\it Step 2:} To conclude the proof, it is enough to show, by Proposition
\ref{p:pfsm-cfsm-minimal-flow}, that the flow $\{\psi_c\}_{c>0}$ is periodic. The idea
can informally be explained as follows. By using (\ref{e:generated}) and
(\ref{e:pfsm-set}), we get that for $c=c(x)\neq 1$,
$$
G(\psi_{c(x)}(x),u)=h(x)G(x,c(x)u+a(x))+j(x)=k(x)G(x,u+b(x))+l(x),
$$
for some $a,b,h\neq 0,j,k\neq 0,l$. Then, for any $t\in\bbR$, $G_t(\Psi(x,u)) =
k(x)G_t(x,u)$, where $G_t$ is defined by (\ref{e:G_t}),
$\Psi(x,u)=(\psi_{c(x)}(x),u-b(x))$ and $k(x)\neq 0$. Since the representation
$\{G_t\}_{t\in\bbR}$ is minimal, $\Psi(x,u)=(x,u)$ and therefore $\psi_{c(x)}(x)=x$ for
$c(x)\neq 1$, showing that the flow $\{\psi_c\}_{c>0}$ is periodic. This argument is not
rigorous because $c$ depends on $x$ and hence the relation (\ref{e:generated}) cannot be
applied directly. The rigorous proof below shows how this technical difficulty can be
overcome.

Consider the set
$$
A = \{(x,c) \in X\times ((0,\infty)\setminus \{1\}) : G(x,cu) = bG(x,u+a) + d\
\mbox{a.e.}\ du
$$
$$
\mbox{for some}\ a=a(x,c),b=b(x,c)\neq 0,d=d(x,c)\in\bbR \}
$$
and let
$$
A_0 = A \cap \{(x,c) \in X\times ((0,\infty)\setminus \{1\}) : G(x,cu) = h
G(\psi_c(x),u+g) + j\ \mbox{a.e.}\ du
$$
$$
\mbox{for some}\ h=h(x,c)\neq 0,g=g(x,c), j = j(x,c)\in\bbR \}.
$$
Since $G$ satisfies (\ref{e:generated}), we have $A=A_0$ a.e.\ $\mu(dx)$ for all $c>0$
and hence, by the Fubini's Theorem (see also Lemma 3.1 in Pipiras and Taqqu
\cite{pipiras:taqqu:2002d}), we have that $A = A_0$ a.e.\ $\mu(dx)\tau(dc)$, where $\tau$
is any $\sigma$-finite measure on $(0,\infty)$. Setting
$$
A_1 = A_0\cap \{(x,c)\in X\times ((0,\infty)\setminus \{1\}) : \psi_c(x) = x\}
$$
we want to show that $A_1=A_0$ a.e.\ $\mu(dx)\tau(dc)$ and to do so, it is enough to
prove that
\begin{equation}\label{e:psi_c(x)=x}
\psi_c(x) = x\quad  \mbox{a.e.}\ \mbox{for}\ (x,c)\in A_0.
\end{equation}

We proceed by contradiction. Suppose that (\ref{e:psi_c(x)=x}) is not true. We can then
find a fixed $c_0\neq 1$ such that $\psi_{c_0}(x)\neq x$ a.e.\ on a set of positive
measure for $(x,c_0)\in A_0$. Define first $\widetilde \psi(x) = \psi_{c_0}(x)$ for
$(x,c_0)\in A_0$ and $\widetilde\psi(x) = x$ for $(x,c_0)\notin A_0$. Then,
\begin{equation}\label{e:contradiction-min1}
    G(\widetilde \psi(x),u+\widetilde a(x)) + \widetilde c(x)= \widetilde b(x) G(x,u)
\end{equation}
a.e.\ $\mu(dx)du$, for some measurable functions $\widetilde a$, $\widetilde b\neq 0$ and
$\widetilde c$. Indeed, relation (\ref{e:contradiction-min1}) is  clearly true for $x$
such that $(x,c_0)\notin A_0$ since $\widetilde \psi(x) = x$. It is also true for
$(x,c_0)\in A_0$ because it follows from the definition of $A_0$ that the relations
$G(x,c_0u) = bG(x,u+a) +d$ and $G(x,c_0u) = h G(\psi_{c_0}(x),u+g)+j$ imply
$G(\psi_{c_0}(x),u+\widetilde a) + \widetilde c = \widetilde bG(x,u)$. Now define
$\Psi(x,u)= (\widetilde \psi(x),u+\widetilde a(x))$. We obtain from
(\ref{e:contradiction-min1}) that, for all $t\in \bbR$,
\begin{equation}\label{e:contradiction-min2}
    G_t(\Psi(x,u)) = h(x) G_t(x,u)\quad \mbox{a.e.}\ \mu(dx)du,
\end{equation}
where $h(x)\neq 0$ and where we used the notation (\ref{e:G_t}). Since $\widetilde \psi$
is nonsingular by construction, the map $\Psi$ is nonsingular as well and, since
$\psi(x)\neq x$ on a set of positive measure $\mu(dx)$, we have $\Psi(x,u)\neq (x,u)$
($\Psi$ is not an identity map) on a set of positive measure $\mu(dx)du$. This
contradicts (\ref{e:G_t-minimal}) and hence the minimality of the representation
$\{G_t\}_{t\in\bbR}$. Hence, $A_1 = A_0$ a.e.\ $\mu(dx)\tau(dc)$ and since $A_0=A$ a.e.\
$\mu(dx)\tau(dc)$ as well, we have
\begin{equation}\label{e:A=A_1}
    A = A_1 \quad \mbox{a.e.}\ \mu(dx)\tau(dc).
\end{equation}
By Lemma 3.3 in Pipiras and Taqqu \cite{pipiras:taqqu:2003cy}, we can choose a
$\mu$-measurable function $c(x)\neq 1$ defined for $x\in\mbox{proj}_X A_1$ such that
$(x,c(x))\in A_1$ and, in particular,
\begin{equation}\label{e:psi-periodic}
    \psi_{c(x)}(x) = x.
\end{equation}
By using (\ref{e:A=A_1}), the definition of $C_P$ and the assumption $C_P = X$
$\mu$-a.e., we obtain that $\mbox{proj}_X A_1 = \mbox{proj}_X A = C_P = X$ $\mu$-a.e.,
that is, (\ref{e:psi-periodic}) holds for $\mu$-a.e.\ $x\in X$. Hence, $X = P$
$\mu$-a.e., showing that the flow $\psi_c$ is periodic.

\bigskip
To prove the converse, suppose that $X_\alpha$ given by (\ref{e:mma}) with a kernel $G$
satisfying (\ref{e:full-support}), is a PFSM. By Proposition
\ref{p:pfsm-cfsm-minimal-flow}, the minimal representation (\ref{e:minimal-mma}) of
$X_\alpha$ is generated by a periodic flow $\{\widetilde \psi_c\}_{c>0}$. Let $\widetilde
P$ be the set of the periodic points of the flow $\{\widetilde \psi_c\}_{c>0}$, and
$\widetilde C_P$ be the PFSM set defined using the representation (\ref{e:minimal-mma}).
Since the flow $\{\widetilde \psi_c\}_{c>0}$ is periodic, $\widetilde P = \widetilde X$
a.e.\ $\widetilde \mu(d\widetilde x)$. Since $\widetilde P\subset \widetilde C_P$ a.e.\
$\widetilde \mu(d\widetilde x)$ by Proposition \ref{p:P-vs-C_P} below, we have
$\widetilde C_P = \widetilde X$ a.e.\ $\widetilde \mu(d\widetilde x)$. In addition, the
following three equalities hold a.e.\ $\mu(dx)$: $C_P = \Phi_1^{-1}(\widetilde C_P)$,
$\Phi_1^{-1}(\widetilde C_P)=\Phi_1^{-1}(\widetilde X)$ and $\Phi_1^{-1}(\widetilde X) =
X$. The first equality follows from (\ref{e:C_P-is-Phi^-1-tile-C_P}), the second holds
because the measures $\mu\circ \Phi_1^{-1}$ and $\widetilde \mu$ are absolutely
continuous by (\ref{e:nonminimal-mu-minimal-mu}) and hence $\widetilde C_X = \widetilde
X$ a.e.\ $\widetilde \mu(d\widetilde x)$ implies $\mu(\Phi_1^{-1}(\widetilde X\setminus
\widetilde C_P)) = 0$. The third equality follows from the definition of $\Phi_1$.
Stringing these equalities together one gets $C_P=X$ a.e.\ $\mu(dx)$. \ \ $\Box$
\end{proof}

\bigskip

The next result describes relations between the PFSM set $C_P$ defined using a kernel
function $G$, and the set of periodic points $P$ of a flow related to the kernel $G$ as
in Definition \ref{d:ss-mma}. The first part of the result was used in the proof of
Theorem \ref{t:identification-pfsm} above.

\begin{proposition}\label{p:P-vs-C_P}
Suppose that a $S\alpha S$, $\alpha\in (0,2)$, self-similar mixed moving average
$X_\alpha$ given by (\ref{e:mma}), is generated by a flow $\{\psi_c\}_{c>0}$. Let $P$ be
the set of periodic points (\ref{e:P}) of the flow $\{\psi_c\}_{c>0}$ and $C_P$ the PFSM
set (\ref{e:pfsm-set}) defined using the kernel $G$ of the representation (\ref{e:mma}).
Then, we have
\begin{equation}\label{e:P-subset-C_P}
 \hfill P\subset C_P\quad \mu\mbox{-a.e.}
\end{equation}
If, moreover, the representation (\ref{e:mma}) is minimal, we have
\begin{equation}\label{e:P=C_P}
 \hfill P = C_P\quad \mu\mbox{-a.e.}
\end{equation}
\end{proposition}

\begin{proof}
We first prove (\ref{e:P-subset-C_P}). Let $\tau(dc)$ denote any $\sigma$-finite measure
on $(0,\infty)$. By the Fubini's theorem (see also Lemma 3.1 in Pipiras and Taqqu
\cite{pipiras:taqqu:2002d}), relation (\ref{e:generated}) implies that a.e.\
$\mu(dx)\tau(dc)$,
$$
G(x,cu) = hG(\psi_c(x),u+g) + j\quad \mbox{a.e.}\ du,
$$
for some $h=h(x,c)\neq 0$, $g=g(x,c)$ and $j=j(x,c)$. Hence, setting
$$
\widetilde P: = \Big\{(x,c)\in X\times ((0,\infty)\setminus \{1\}): \psi_c(x) = x\Big\},
$$
we have a.e.\ $\mu(dx)\tau(dc)$,
$$
\widetilde P  = \widetilde P \bigcap \Big\{(x,c): G(x,cu) = hG(\psi_c(x),u+g)+j\
\mbox{a.e.}\ du\ \mbox{for some}\ h\neq 0,g,j \Big\}
$$
\begin{equation}\label{e:tildeP=tildeP-int}
= \widetilde P \bigcap \Big\{(x,c): G(x,cu) = hG(x,u+g)+j\ \mbox{a.e.}\ du\ \mbox{for
some}\ h\neq 0,g,j \Big\}.
\end{equation}
Since $P = \mbox{proj}_X\widetilde P$, relation (\ref{e:tildeP=tildeP-int}) implies that
a.e.\ $x\in P$ belongs to the set
$$
\mbox{proj}_X \Big( \Big\{(x,c): G(x,cu) = hG(x,u+g)+j\ \mbox{a.e.}\ du\ \mbox{for some}\
h\neq 0,g,j \Big\}\Big),
$$
that is, there is $c=c(x)\neq 1$ such that
$$
G(x,cu) = hG(x,u+g) + j\quad \mbox{a.e.}\ du
$$
for some $h\neq 0$,$g$,$j$. This shows that $P\subset C_P$ a.e.\ $\mu(dx)$.

\smallskip
To prove (\ref{e:P=C_P}), suppose that the representation (\ref{e:mma}) is minimal. It is
enough to show that $C_P\subset P$ $\mu$-a.e. Let $\{G_t\vert_{C_P}\}$ be the kernel
$G_t$ of (\ref{e:mma}) restricted to the set $C_P\times \bbR$. By Lemma
\ref{l:C_P-invariant} below, the set $C_P$ is a.e.\ invariant under the flow
$\{\psi_c\}$. Then, $\{G_t\vert_{C_P}\}$ is a representation of a self-similar mixed
moving average. Since $\{G_t\}$ is minimal, so is the representation
$\{G_t\vert_{C_P}\}$. It is obviously generated by the flow $\psi\vert_{C_P}$, the
restriction of the flow $\psi$ to the set $C_P$. By arguing as in Step 2 of the proof of
Theorem \ref{t:identification-pfsm}, we therefore obtain that for a.e.\ $x\in C_P$,
$\psi_{c(x)}(x)=x$ for some $c(x)\neq 1$. This shows that $C_P\subset P$ a.e.\ $\mu(dx)$.
\ \ $\Box$
\end{proof}

\medskip

The following lemma was used in the proof of Proposition \ref{p:P-vs-C_P} above.

\begin{lemma}\label{l:C_P-invariant}
If a $S\alpha S$, $\alpha\in (0,2)$, self-similar mixed moving average $X_\alpha$ given
by a representation (\ref{e:mma}) is generated by a flow $\{\psi_c\}_{c>0}$, and $C_P$ is
the PFSM set defined by (\ref{e:pfsm-set}), then $C_P$ is a.e.\ invariant under the flow
$\{\psi_c\}_{c>0}$, that is, $\mu(C_P\triangle \psi_c^{-1}(C_P)) = 0$ for all $c>0$.
\end{lemma}

\begin{proof}
Since $\{\psi_c\}_{c>0}$ satisfies the group property (\ref{e:flow0}), it is enough to
show that $C_P \subset \psi_{r}^{-1}(C_P)$ $\mu$-a.e.\ for any fixed $r>0$. By
(\ref{e:generated}), we have for any $c>0$,
\begin{equation}\label{e:for-C_P-invariant1}
G(\psi_{r}(x),cu+a(x)) = b(c,x) G(x,cru) + j(c,x)\quad \mbox{a.e.}\ \mu(dx)du,
\end{equation}
for some $a,b\neq 0,j$ (these depend on $r$ but since $r$ is fixed we do not indicate
their dependence on $r$). By using Lemma 4.2 in Pipiras and Taqqu
\cite{pipiras:taqqu:2003cy} and arguing as in Step 2 of the proof of Theorem
\ref{t:identification-pfsm}, we can choose a function $c(x)\neq 1$ such that, for a.e.\
$x\in C_P$,
\begin{equation}\label{e:for-C_P-invariant2}
G(x,c(x)ru) = b(x) G(x,ru+a(x)) + j(x)\quad \mbox{a.e.}\ du,
\end{equation}
for some $a,b\neq 0,j$, and such that the relation (\ref{e:for-C_P-invariant1}) holds
with $c$ replaced by $c(x)$. By substituting (\ref{e:for-C_P-invariant2}) into
(\ref{e:for-C_P-invariant1}) with $c=c(x)$ and then making a change of variables in $u$,
we obtain that, for a.e.\ $x\in C_P$,
$$
G(\psi_{r}(x),c(x)u+d(x)) = h(x) G(x,ru) + l(x)\quad \mbox{a.e.}\ du,
$$
for some $d,h\neq 0,l$. Then, by using (\ref{e:generated}) and making a change of
variables in $u$, we get that, for a.e.\ $x\in C_P$,
$$
G(\psi_{r}(x),c(x)u) = k(x) G(\psi_r(x),u+p(x)) + q(x)\quad \mbox{a.e.}\ du,
$$
for some $k\neq 0,p,q$. Hence, for a.e.\ $x\in C_P$, $\psi_{r}(x)\in C_P$ or $x\in
\psi_{r}^{-1}(C_P)$, showing that $C_P\subset \psi_{r}^{-1}(C_P)$ $\mu$-a.e. \ \ $\Box$
\end{proof}

\medskip

We now provide an example of a PFSM. Further examples of PFSMs can be found in Pipiras
and Taqqu \cite{pipiras:taqqu:2003ca}.

\begin{example}\label{ex:periodic}
Let $\alpha\in (0,2)$, $H\in(0,1)$ and $\kappa = H - 1/\alpha<0$. By Section 8 of Pipiras
and Taqqu \cite{pipiras:taqqu:2002s}, the mixed moving average process
\begin{equation}\label{e:ex-periodic}
\int_0^1\int_\bbR \Big((t+u)_+^\kappa 1_{[0,1/2)}(\{x+\ln|t+u|\}) -  u_+^\kappa
1_{[0,1/2)}(\{x+\ln|u|\})\Big) M_\alpha(dx,du),
\end{equation}
where $M_\alpha$ has the control measure $dxdu$ and $\{u\}$ stands for the fractional
part of $u\in\bbR$, is well-defined and self-similar. It has the representation
(\ref{e:mma}) with $X=[0,1)$ and
$$
G(x,u) = u_+^\kappa 1_{[0,1/2)}(\{x+\ln|u|\}),\ x\in[0,1),u\in\bbR.
$$
Since $G(x,eu) = e^\kappa G(x,u)$ for all $x\in[0,1),u\in\bbR$, we deduce that $X=C_P$
for the process (\ref{e:ex-periodic}). Hence, by Theorem \ref{t:identification-pfsm}, the
process (\ref{e:ex-periodic}) is a PFSM when $\alpha \in (1,2)$.
\end{example}


\section{Identification of cyclic fractional stable motions: the non-minimal
case}\label{s:cfsm-nonminimal}

We focused so far on periodic fractional stable motions. By using the set $C_P$ we were
able to identify them without requiring the representation to be minimal. We now want to
do the same thing for cyclic fractional stable motions by introducing a corresponding set
$C_L$. To do so, observe that

\begin{lemma}\label{l:cfsm=pfsm+no-mlfsm}
A CFSM is a PFSM without a mixed LFSM component.
\end{lemma}

\begin{proof}
This follows from (\ref{e:into-four-minimal}) and the fact that $X_\alpha^F$ is a mixed
LFSM (see (\ref{e:mlfsm})). \ \ $\Box$
\end{proof}

\medskip
We showed in Pipiras and Taqqu \cite{pipiras:taqqu:2002s} that a mixed LFSM can be
identified through the {\it mixed LFSM set}
\begin{eqnarray}
  C_F &=& \Big\{ x\in X: G(x,u) = d(u+f)_+^\kappa + h(u+f)_-^\kappa + g\ \mbox{a.e.}\ du \nonumber \\
   & & \hspace{0.3in} \mbox{for some reals}\ d = d(x),f=f(x),g=g(x),h=h(x) \Big\}   \label{e:C_F-not0}
\end{eqnarray}
when $\kappa\neq 0$, and
\begin{eqnarray}
  C_F &=& \Big\{ x\in X: G(x,u) = d\ln|u+f| + h 1_{(0,\infty)}(u+f) + g\ \mbox{a.e.}\ du \nonumber \\
   & & \hspace{0.3in} \mbox{for some reals}\ d = d(x),f=f(x),g=g(x),h=h(x) \Big\}   \label{e:C_F-0}
\end{eqnarray}
when $\kappa= 0$. The following lemma shows that this set is a subset of $C_P$.

\begin{lemma}\label{l:C_F-subset-C_P}
We have
\begin{equation}\label{e:C_F-subset-C_P}
C_F\subset C_P.
\end{equation}
\end{lemma}

\begin{proof}
Suppose that $\kappa\neq 0$. If $x\in C_F$, then $G(x,u) = d(u+f)_+^\kappa +
h(u+f)_-^\kappa + g$ for some reals $d,f,g,h$ and hence
$$
G(x,cu) = c^{\kappa} \Big( d(u+c^{-1}f)_+^\kappa + h(u+c^{-1}f)_-^\kappa + g\Big) +
(1-c^\kappa)g
$$
\begin{equation}\label{e:mixed-lfsm-Gc}
= c^{\kappa}G(x,u+c^{-1}f) + (1-c^\kappa)g
\end{equation}
for arbitrary $c$. This shows that $x\in C_P$ and hence that (\ref{e:C_F-subset-C_P})
holds. The proof in the case $\kappa= 0$ is similar.  \ \ $\Box$
\end{proof}

\medskip Since a CFSM is a PFSM without a mixed LFSM component, we expect that a CFSM can be
identified through the set $C_L = C_P\setminus C_F$. We will show that this is indeed the
case.

\begin{definition}\label{d:cfsm-set}
A {\it cyclic fractional stable motion set} ({\it CFSM set}, in short) of a self-similar
mixed moving average $X_\alpha$ given by (\ref{e:mma}) is defined by
\begin{equation}\label{e:cfsm-set}
    C_L := C_P\setminus C_F,
\end{equation}
where $C_P$ is the PFSM set defined by (\ref{e:pfsm-set}) and $C_F$ is the mixed LFSM set
defined by (\ref{e:C_F-not0})--(\ref{e:C_F-0}).
\end{definition}

The following result shows that a CFSM can indeed be identified through the CFSM set.

\begin{theorem}\label{t:identification-cfsm}
A $S\alpha S$, $\alpha\in (1,2)$, self-similar mixed moving average $X_\alpha$ given by
(\ref{e:mma}) with $G$ satisfying (\ref{e:full-support}), is a CFSM if and only if $C_L =
X$ $\mu$-{\rm a.e.}, where $C_L$ is the CFSM set defined in (\ref{e:cfsm-set}).
\end{theorem}

\begin{proof}
If $X_\alpha$ is a CFSM, then it is also a PFSM and hence, by Theorem
\ref{t:identification-pfsm}, $C_P=X$ $\mu$-a.e. By (\ref{e:cfsm-set}), $C_P = C_L+C_F$.
Since $X_\alpha$ does not have a mixed LFSM component (Lemma \ref{l:cfsm=pfsm+no-mlfsm}
above), Propositions 7.1 and 7.2 in Pipiras and Taqqu \cite{pipiras:taqqu:2002s} imply
that $C_F=\emptyset$ $\mu$-a.e. Hence, $C_L=X$ $\mu$-a.e. Conversely, if $C_L=X$
$\mu$-a.e., then $C_P=X$ $\mu$-a.e.\ and hence $X_\alpha$ is a PFSM. But $C_L=X$
$\mu$-a.e.\ implies $C_F=\emptyset$ $\mu$-a.e., that is, $X_\alpha$ does not have a mixed
LFSM component. The PFSM $X_\alpha$ is therefore a CFSM. \ \ $\Box$
\end{proof}

\medskip Observe that the mixed LFSM set $C_F$ in (\ref{e:C_F-not0})--(\ref{e:C_F-0}) is expressed in a
different way from the PFSM set (\ref{e:pfsm-set}). It can, however, be expressed in a
similar way.

\begin{proposition}\label{p:C_F-alternative}
Let $\alpha \in (1,2)$. We have
\begin{eqnarray}
C_F & = & \Big\{ x\in X: \exists\ c_n = c_n(x)\to 1\ (c_n\neq 1): G(x,c_nu)
= b_n\, G(x,u + a_n) + d_n\ \mbox{a.e.}\ du\nonumber \\
& &\hspace{.4in} \mbox{for some}\ a_n = a_n(c_n,x), b_n=b_n(c_n,x)\neq
0,d_n=d_n(c,x)\in\bbR \Big\}, \quad \mu\mbox{-a.e.}\quad \label{e:C_F-alternative}
\end{eqnarray}
\end{proposition}

\begin{proof}
Consider the case $\kappa\neq 0$. Denote the set on the right-hand side of
(\ref{e:C_F-alternative}) by $C_F^0$. If $x\in C_F$, then for any $c\neq 1$, $G(x,cu) =
c^\kappa G(x,u+c^{-1}f)+(1-c^\kappa)g$ (see (\ref{e:mixed-lfsm-Gc})) and hence $x\in
C_F^0$ with any $c_n\to 1$ ($c_n\neq 1$). This shows that $C_F\subset C_F^0$ in the case
$\kappa\neq 0$. The proof in the case $\kappa  =0$ is similar.

To show that $C_F^0\subset C_F$ $\mu$-a.e., we adapt the proof of Proposition 5.1 in
Pipiras and Taqqu \cite{pipiras:taqqu:2003cy}. Let $\widetilde G:\widetilde X\times\bbR
\mapsto \bbR$ be the kernel function of a minimal representation of the process
$X_\alpha$, and $\widetilde C_F$ and $\widetilde C_F^0$ be the sets defined in the same
way as $C_F$ and $C_F^0$ by using the kernel function $\widetilde G$. One can show as in
the proof of Theorem \ref{t:identification-pfsm} (see (\ref{e:C_P-is-Phi^-1-tile-C_P}))
that $C_F^0 = \Phi_1^{-1}(\widetilde C_F^0)$ $\mu$-a.e., where $\Phi_1$ is the map
appearing in (\ref{e:nonminimal-minimal}) and (\ref{e:nonminimal-mu-minimal-mu}). As
shown in the proof of Proposition 7.1 of Pipiras and Taqqu \cite{pipiras:taqqu:2002s},
$C_F = \Phi^{-1}_1(\widetilde C_F)$. By using (\ref{e:nonminimal-mu-minimal-mu}), it is
then enough to show that $\widetilde C_F^0\subset \widetilde C_F$ $\widetilde \mu$-a.e.,
or equivalently, $C_F^0\subset C_F$ $\mu$-a.e.\ but where $C_F^0$ and $C_F$ are defined
by using the kernel function $G$ corresponding to a minimal representation of $X_\alpha$.

If the process $X_\alpha$ is given by a minimal representation involving a kernel $G$,
then it is generated by a flow $\{\psi_c\}_{c>0}$ and related functionals (Theorem 4.1 in
Pipiras and Taqqu \cite{pipiras:taqqu:2002d}). By Lemma \ref{l:C_F^0-invariant} below,
the set $C_F^0$ is a.e.\ invariant under the flow $\{\psi_c\}_{c>0}$. Then, the process
\begin{equation}\label{e:X_alpha-over-C_F^0}
\int_{C_F^0}\int_\bbR \Big(G(x,t+u) - G(x,u) \Big)M_\alpha(dx,du)
\end{equation}
is a self-similar mixed moving average, the representation (\ref{e:X_alpha-over-C_F^0})
is minimal and the process (\ref{e:X_alpha-over-C_F^0}) is generated by the flow
$\{\psi_c\}_{c>0}$ restricted to the set $C_F^0$. Arguing as in the proof of Theorem
\ref{t:identification-pfsm}, one can show that, for a.e.\ $x\in C_F^0$,
\begin{equation}\label{e:psi-c_n-to-1}
\psi_{c_n(x)}(x) = x\quad \mbox{for}\ \ c_n(x)\to 1\ (c_n(x)\neq 1).
\end{equation}
Relation (\ref{e:psi-c_n-to-1}) cannot hold for points which are not fixed. This follows
by using the so-called ``special representation'' of a flow as in the end of the proof of
Proposition 5.1, Pipiras and Taqqu \cite{pipiras:taqqu:2003cy} (see relation (5.6) of
that paper). Hence, for a.e.\ $x\in C_F^0$, $\psi_c(x)=x$ for all $c>0$. Since $C_F = F =
\{x:\psi_c(x) = x\ \mbox{for all} \ c>0\}$ $\mu$-a.e.\ by Theorem 10.1 in Pipiras and
Taqqu \cite{pipiras:taqqu:2002s}, we obtain that $C_F^0\subset C_F$ $\mu$-a.e. \ \ $\Box$
\end{proof}

\medskip
The following lemma was used in the proof of Proposition \ref{p:C_F-alternative}. An
a.e.\ invariant set is defined in Lemma \ref{l:C_P-invariant} above.

\begin{lemma}\label{l:C_F^0-invariant}
If a $S\alpha S$, $\alpha\in (0,2)$, self-similar mixed moving average $X_\alpha$ given
by a representation (\ref{e:mma}) is generated by a flow, and $C_F^0$ denotes the
right-hand side of (\ref{e:C_F-alternative}), then $C_F^0$ is a.e.\ invariant under the
flow.
\end{lemma}

\begin{proof}
Since the proof of this result is very similar to that of Lemma \ref{l:C_P-invariant}, we
only outline it. Proceeding as in the proof of Lemma \ref{l:C_P-invariant}, we can choose
functions $c_n(x)\to 1$ ($c_n(x)\neq 1$) such that, for a.e.\ $x\in C_F^0$, the relation
(\ref{e:for-C_P-invariant2}) holds with $c(x)$ replaced by $c_n(x)$ (and with
$a_n,b_n\neq 0,j_n$ replacing $a,b\neq 0,j$) and the relation
(\ref{e:for-C_P-invariant1}) holds with $c$ replaced by $c_n(x)$. The conclusion follows
as in the proof of Lemma \ref{l:C_P-invariant}. \ \ $\Box$
\end{proof}

\medskip
The new formulation (\ref{e:C_F-alternative}) of $C_F$ yields the following
characterization of $C_L =C_P\setminus C_F$:

\begin{corollary}\label{c:cfsm-set-alternative}
We have
\begin{eqnarray}
C_L & = & \Big\{ x\in X: \exists\ c_0 = c_0(x) \neq 1, \nexists\ c_n =
c_n(x) \to 1\ (c_n\neq 1): \nonumber \\
& & \hspace{.4in} G(x,c_n u) = b_n\, G(x,u + a_n) + d_n\  \mbox{a.e.}\ du,\ n=0,1,2,\ldots\nonumber \\
& &\hspace{.5in} \mbox{for some}\ a_n = a_n(c_n,x), b_n=b_n(c_n,x)\neq
0,d_n=d_n(c,x)\in\bbR \Big\}, \quad \mu\mbox{-a.e.}\quad \label{e:cfsm-alternative}
\end{eqnarray}
\end{corollary}

The next result is analogous to the second part of Proposition \ref{p:P-vs-C_P}.

\begin{proposition}\label{p:L=C_L}
Suppose that a $S\alpha S$, $\alpha\in (0,2)$, self-similar mixed moving average
$X_\alpha$ given by a minimal representation (\ref{e:mma}), is generated by a flow
$\{\psi_c\}_{c>0}$. Then,
\begin{equation}\label{e:L=C_L}
 \hfill L = C_L\quad \mu\mbox{-a.e.},
\end{equation}
where $L$ is the set of cyclic points (\ref{e:L}) of the flow $\{\psi_c\}_{c>0}$ and
$C_L$ is the CFSM set (\ref{e:cfsm-set}) defined using the kernel of a minimal
representation (\ref{e:mma}).
\end{proposition}

\begin{proof}
By Proposition \ref{p:P-vs-C_P} above and Theorem 10.1 in Pipiras and Taqqu
\cite{pipiras:taqqu:2002s}, we have $P=C_P$ $\mu$-a.e.\ and $F=C_F$ $\mu$-a.e., where $P$
and $F$ are the sets of the periodic and fixed points of the flow $\{\psi_c\}_{c>0}$, and
$C_P$ and $C_F$ are the PFSM and the mixed LFSM sets. The equality (\ref{e:L=C_L})
follows since $L=P\setminus F$ and $C_L=C_P\setminus C_F$.  \ \ $\Box$
\end{proof}

\medskip
The PFSM considered in Example \ref{ex:periodic} above is also a CFSM.

\begin{example}
The self-similar mixed moving average (\ref{e:ex-periodic}) considered in Example
\ref{ex:periodic} above is a CFSM because it is a PFSM and, as can be seen by using
(\ref{e:C_F-0}), $C_F = \emptyset$.
\end{example}


\section{Refined decomposition of self-similar mixed moving averages}
\label{s:refined-decomposition}

Suppose that $X_\alpha$ is a $S\alpha S$, $\alpha\in (1,2)$, self-similar mixed moving
average. By using its minimal representation, we showed in Section \ref{s:pcfsm-minimal}
that $X_\alpha$ admits a decomposition (\ref{e:into-four-minimal}) which is unique in
distribution and has independent components. We show here that the components of the
decomposition (\ref{e:into-four-minimal}) can be expressed in terms of a possibly
nonminimal representation (\ref{e:mma}) of the process $X_\alpha$.

Let $G$ be the kernel function of a possibly nonminimal representation (\ref{e:mma}) of
the process $X_\alpha$. With the notation (\ref{e:G_t}), let
\begin{eqnarray}
  D  &=& \Big\{x\in X : \int_0^\infty dc \int_\bbR du\ c^{-H\alpha} |G_c(x,cu)|^\alpha < \infty   \Big\}, \label{e:D} \\
  C  &=& \Big\{x\in X : \int_0^\infty dc \int_\bbR du\ c^{-H\alpha} |G_c(x,cu)|^\alpha = \infty   \Big\}. \label{e:C}
\end{eqnarray}
Recall also the definitions (\ref{e:C_F-not0})--(\ref{e:C_F-0}), (\ref{e:pfsm-set}) and
(\ref{e:cfsm-set}) of the mixed LFSM, PFSM and CFSM sets defined by using the kernel
function $G$.

\begin{theorem}\label{t:into-four}
Let $X_\alpha$ be a $S\alpha S$, $\alpha\in (1,2)$, self-similar mixed moving average
given by a possibly nonminimal representation (\ref{e:mma}). Suppose that
$$
X_\alpha^D,\ X_\alpha^F,\ X_\alpha^L,\ X_\alpha^{C\setminus P}
$$
are the four independent components in the unique decomposition
(\ref{e:into-four-minimal}) of the process $X_\alpha$ obtained by using its minimal
representation. Then,
\begin{eqnarray}
  X_\alpha^D(t) &\stackrel{d}{=}& \int_D \int_\bbR G_t(x,u) M_\alpha(dx,du), \label{e:int-D=D} \\
  X_\alpha^F(t) &\stackrel{d}{=}& \int_{C_F} \int_\bbR G_t(x,u) M_\alpha(dx,du), \label{e:int-F=C_F} \\
  X_\alpha^L(t) &\stackrel{d}{=}& \int_{C_L} \int_\bbR G_t(x,u) M_\alpha(dx,du), \label{e:int_L=C_L} \\
  X_\alpha^{C\setminus P}(t) &\stackrel{d}{=}& \int_{C\setminus C_P}\int_\bbR G_t(x,u)
  M_\alpha(dx,du),\label{e:int_C-P=C-C_P}
\end{eqnarray}
where $\stackrel{d}{=}$ stands for the equality in the sense of the finite-dimensional
distributions and the sets $D$, $C$, $C_F$, $C_P$ and $C_L$ are defined by (\ref{e:D}),
(\ref{e:C}), (\ref{e:C_F-not0})--(\ref{e:C_F-0}), (\ref{e:pfsm-set}) and
(\ref{e:cfsm-set}), respectively.
\end{theorem}

\begin{proof}
The equalities (\ref{e:int-D=D}) and (\ref{e:int-F=C_F}) follow from Theorem 5.5 in
Pipiras and Taqqu \cite{pipiras:taqqu:2002d} and Corollary 9.1 in Pipiras and Taqqu
\cite{pipiras:taqqu:2002s}, respectively. Consider now the equality (\ref{e:int_L=C_L}).
Let $\widetilde G$ be the kernel of a minimal representation (\ref{e:minimal-mma}) of the
process $X_\alpha$, and let also $\widetilde C_F$, $\widetilde C_P$ and $\widetilde C_L$
be the sets defined by (\ref{e:C_F-not0})--(\ref{e:C_F-0}), $(\ref{e:pfsm-set})$ and
(\ref{e:cfsm-set}), respectively, using the kernel function $\widetilde G$. Since $C_P =
\Phi_1^{-1}(\widetilde C_P)$ $\mu$-a.e.\ by (\ref{e:C_P-is-Phi^-1-tile-C_P}) and $C_F =
\Phi_1^{-1}(\widetilde C_F)$ $\mu$-a.e.\ as shown in the proof of Proposition 7.1 in
Pipiras and Taqqu \cite{pipiras:taqqu:2002s}, we obtain that $C_L = C_P\setminus C_F =
\Phi_1^{-1}(\widetilde C_P\setminus \widetilde C_F) =  \Phi_1^{-1}(\widetilde C_L)$
$\mu$-a.e. Then, by using (\ref{e:nonminimal-minimal}),
(\ref{e:nonminimal-mu-minimal-mu}) and a change of variables as at the end of the proof
of Proposition 7.1 in Pipiras and Taqqu \cite{pipiras:taqqu:2002s}, we get that
$$
\int_{C_L} \int_\bbR G_t(x,u) M_\alpha(dx,du)  \stackrel{d}{=}
  \int_{\widetilde C_L}\int_\bbR \widetilde G_t(\widetilde x,u) \widetilde
  M_\alpha(d\widetilde x,du).
$$
Since $\widetilde G$ is a kernel of a minimal representation, it is related to a flow in
the sense of Definition \ref{d:ss-mma}. Let $\widetilde L$ be the set of the cyclic
points of the flow corresponding to the kernel $\widetilde G$. Since $\widetilde L =
\widetilde C_L$ $\mu$-a.e.\ by Proposition \ref{p:L=C_L}, we get that
\begin{equation}\label{e:int_C_L=tilde-C_L}
\int_{C_L} \int_\bbR G_t(x,u) M_\alpha(dx,du)  \stackrel{d}{=} \int_{\widetilde L}
\int_\bbR \widetilde G_t(\widetilde x,u) \widetilde
  M_\alpha(d\widetilde x,du).
\end{equation}
The process on the right-hand side of (\ref{e:int_C_L=tilde-C_L}) has the distribution of
$X_\alpha^L$ by the definition of $X_\alpha^L$ and the uniqueness result in Theorem
\ref{t:unique-minimal}.

To show the equality (\ref{e:int_C-P=C-C_P}), observe that by Lemma
\ref{l:C_F-subset-C_P} and Lemma \ref{l:C_P-subset-C} below, we have $C_F\subset
C_P\subset C$. Since $C_P = C_F + C_L$, the sets $C_F$, $C_L$ and $C\setminus C_P$ are
disjoint, and $C_F+C_L+C\setminus C_P = C$. Hence, the processes on the right-hand side
of (\ref{e:int-D=D})--(\ref{e:int_C-P=C-C_P}) are independent. Since the processes on the
left-hand side of (\ref{e:int-D=D})--(\ref{e:int_C-P=C-C_P}) are also independent, since
the sum of the processes on the left-hand side of
(\ref{e:int-D=D})--(\ref{e:int_C-P=C-C_P}) has the same distribution as the sum of the
processes on the right-hand side of (\ref{e:int-D=D})--(\ref{e:int_C-P=C-C_P}), and since
we already showed that the equalities (\ref{e:int-D=D})--(\ref{e:int_L=C_L}) hold, we
conclude that the equality (\ref{e:int_C-P=C-C_P}) holds as well. \ \ $\Box$
\end{proof}

\medskip
The following lemma was used in the proof of Theorem \ref{t:into-four} above.

\begin{lemma}\label{l:C_P-subset-C}
We have
\begin{equation}\label{e:C_P-subset-C}
C_P\subset C,
\end{equation}
where $C_P$ is the PFSM set (\ref{e:pfsm-set}) and $C$ is defined by (\ref{e:C}).
\end{lemma}

\begin{proof}
If $x\in C_P$, then by (\ref{e:G_t}) and (\ref{e:pfsm-set}),
$$
G_{rc}(x,rcu) = G(x,rc(1+u)) - G(x,rcu)
$$
$$
= b(G(x,c(1+u)+a) - G(x,cu+a)) = bG_c(x,cu+a)\quad \mbox{a.e.}\ du,
$$
for any $c>0$ and some $r=r(x)\neq 1$, $b=b(x)\neq 0$ and $a=a(x)$. Suppose without loss
of generality that $r=r(x)>1$. Then, by making changes of variables $c$ to $rc$ and $u$
to $u-c^{-1}a$, we obtain that, for any $n\in\bbZ$,
$$
\int_{r^n}^{r^{n+1}} dc \int_\bbR du\, c^{-H\alpha}|G_c(x,cu)|^\alpha = r^{1-H\alpha}
|b|^\alpha \int_{r^{n-1}}^{r^n} dc \int_\bbR du\, c^{-H\alpha}|G_c(x,cu)|^\alpha
$$
and hence
$$
\int_{r^n}^{r^{n+1}} dc \int_\bbR du\, c^{-H\alpha}|G_c(x,cu)|^\alpha = r^{(1-H\alpha)n}
|b|^{\alpha n} \int_{1}^{r} dc \int_\bbR du\, c^{-H\alpha}|G_c(x,cu)|^\alpha.
$$
This yields that, for $x\in C_P$,
$$
\int_0^\infty dc\int_\bbR du\, c^{-H\alpha} |G_c(x,cu)|^\alpha = \sum_{n=-\infty}^\infty
\int_{r(x)^n}^{r(x)^{n+1}} dc \int_\bbR du\, c^{-H\alpha}|G_c(x,cu)|^\alpha
$$
$$
= \int_1^{r(x)} dc \int_\bbR\, c^{-H\alpha}|G_c(x,cu)|^\alpha du \sum_{n=-\infty}^\infty
r(x)^{(1-H\alpha)n} |b(x)|^{n\alpha} =\infty,
$$
since $\sum_{n=-\infty}^0 r^{(1-H\alpha)n} |b|^{n\alpha} + \sum_{n=1}^{\infty}
r^{(1-H\alpha)n} |b|^{n\alpha}=\infty$, which shows that $x\in C$. \ \ $\Box$
\end{proof}

\medskip
The following theorem is essentially a reformulation of Theorem \ref{t:into-four} and
some other previous results. It provides a decomposition of self-similar mixed moving
averages which is more refined than those established in Pipiras and Taqqu
\cite{pipiras:taqqu:2002d,pipiras:taqqu:2002s}. As in Section \ref{s:pcfsm-minimal}, we
will say that a decomposition of a process $X_\alpha$ obtained from its representation
(\ref{e:mma}) is unique in distribution if the distribution of its components does not
depend on the representation (\ref{e:mma}). We will also say that a process does not have
a PFSM component if it cannot be expressed as the sum of two independent processes where
one process is a PFSM.

\begin{theorem}\label{t:into-four-final}
Let $X_\alpha$ be a $S\alpha S$, $\alpha\in (1,2)$, self-similar mixed moving average
given by a possibly nonminimal representation (\ref{e:mma}). Then, the process $X_\alpha$
can be decomposed uniquely in distribution into four independent processes
\begin{equation}\label{e:into-four-final}
    X_\alpha\stackrel{d}{=} X_\alpha^D + X_\alpha^F + X_\alpha^L + X_\alpha^{C\setminus
    P},
\end{equation}
where
\begin{eqnarray}
  X_\alpha^D(t) &=& \int_D \int_\bbR G_t(x,u) M_\alpha(dx,du), \label{e:X_alpha^D} \\
  X_\alpha^F(t) &=& \int_{C_F} \int_\bbR G_t(x,u) M_\alpha(dx,du), \label{e:X_alpha^F} \\
  X_\alpha^L(t) &=& \int_{C_L} \int_\bbR G_t(x,u) M_\alpha(dx,du), \label{e:X_alpha^L} \\
  X_\alpha^{C\setminus P}(t) &=& \int_{C\setminus C_P}\int_\bbR G_t(x,u)
  M_\alpha(dx,du),\label{e:X_alpha^C-P}
\end{eqnarray}
and the sets $D$, $C$, $C_F$, $C_P$ and $C_L$ are defined by (\ref{e:D}), (\ref{e:C}),
(\ref{e:C_F-not0})--(\ref{e:C_F-0}), (\ref{e:pfsm-set}) and (\ref{e:cfsm-set}),
respectively. Here:

\smallskip
$(i)$ The process $X_\alpha^D$ has the canonical representation given in Theorem 4.1 of
Pipiras and Taqqu \cite{pipiras:taqqu:2002s}, and is generated by a dissipative flow.

\smallskip
$(ii)$ The process $X_\alpha^F$ is a mixed LFSM and has the representation
(\ref{e:mlfsm}).

\smallskip
$(iii)$ The process $X_\alpha^L$ is a CFSM, and the sum $X_\alpha^P=X_\alpha^F +
X_\alpha^L$ is a PFSM.

\smallskip
$(iv)$ The process $X_\alpha^{C\setminus P}$ is a self-similar mixed moving average
without a PFSM component.

\medskip
If the process $X_\alpha$ is generated by a flow $\{\psi_c\}_{c>0}$ then the sets $D$ and
$C$ are identical (a.e.) to the dissipative and the conservative parts of the flow
$\{\psi_c\}_{c>0}$.

If, in addition, the representation of the process $X_\alpha$ is minimal, then the sets
$C_P$, $C_F$ and $C_L$ are the sets of the periodic, fixed and cyclic points of the flow
$\{\psi_c\}_{c>0}$, respectively.
\end{theorem}

\smallskip
\noindent {\bf Remark.} It is important to distinguish
(\ref{e:X_alpha^D})--(\ref{e:X_alpha^C-P}) from
(\ref{e:int-D=D})--(\ref{e:int_C-P=C-C_P}). Because of the relations
(\ref{e:int-F=C_F})--(\ref{e:int_C-P=C-C_P}), the processes $X_\alpha^F$, $X_\alpha^L$
and $X_\alpha^{C\setminus P}$ defined through (\ref{e:X_alpha^F})--(\ref{e:X_alpha^C-P})
are equal in finite-dimensional distributions with the corresponding processes
$X_\alpha^F$, $X_\alpha^L$ and $X_\alpha^{C\setminus P}$ defined through
(\ref{e:X_alpha^S}). They are not identical to them because we are integrating here with
respect to the sets $C_F$, $C_L$ and $C\setminus C_P$ which are defined in terms of the
kernel $G$ whereas in the integration in (\ref{e:X_alpha^S}), one is integrating with
respect to the sets $F$, $L$ and $C\setminus P$ which are defined in terms of the flow
$\{\psi_c\}_{c>0}$. We use the same notation for convenience. The abuse is small because
one has equality in distribution and because $C_F=F$, $C_L=L$ and $C_P=P$ when working
with minimal representations.

In the case of the process $X_\alpha^D$ defined through (\ref{e:int-D=D}), the notation
is consistent because $D$, defined by (\ref{e:D}) in terms of the kernel function $G$, is
equal to the set of dissipative points of the flow $\{\psi_c\}_{c>0}$ for arbitrary, not
necessarily minimal, representations (see Corollary 5.2 in Pipiras and Taqqu
\cite{pipiras:taqqu:2002d}).

\medskip
\begin{proof}
The uniqueness of the decomposition (\ref{e:into-four-final}) into four independent
component follows by using Theorem \ref{t:into-four} and the uniqueness result in Theorem
\ref{t:unique-minimal}. Parts $(i)$ and $(ii)$ follow from Theorem 9.1 in Pipiras and
Taqqu \cite{pipiras:taqqu:2002d}. Part $(iii)$ is a consequence of the equalities
(\ref{e:int-F=C_F}) and (\ref{e:int_L=C_L}) in Theorem \ref{t:into-four} and Definition
\ref{d:pfsm-cfsm-minimal}. To show that the process $X_\alpha^{C\setminus P}$ does not
have a PFSM component, we argue by contradiction. Suppose on the contrary that
$X_\alpha^{C\setminus P}$ has a PFSM component, that is,
$$
X_\alpha^{C\setminus P} \stackrel{d}{=} V + W,
$$
where $V$ and $W$ are independent, and $W$ is a PFSM. Let $G^{C\setminus P}:(C\setminus
P)\times \bbR\mapsto \bbR$ and $F:Y\times \bbR\mapsto \bbR$ be the kernel functions in
the representation of $X_\alpha^{C\setminus P}$ and $W$, respectively, where the integral
representation  of $W$ is equipped with the control measure $\sigma(dy)du$. By using
Theorem 5.2 in Pipiras and Taqqu \cite{pipiras:taqqu:2002d}, there are functions
$\Phi_1:Y\mapsto C\setminus C_P$, $h:Y\mapsto \bbR\setminus\{0\}$ and
$\Phi_2,\Phi_3:Y\mapsto \bbR$ such that
\begin{equation}\label{e:fourth-no-periodic}
F(y,u) = h(y) G^{C\setminus P}(\Phi_1(y),u+\Phi_2(y)) + \Phi_3(y),\quad \mbox{a.e.}\
\sigma(dy)du
\end{equation}
or
\begin{equation}\label{e:fourth-no-periodic2}
G^{C\setminus P}(\Phi_1(y),u) = (h(y))^{-1} F(y,u-\Phi_2(y)) - (h(y))^{-1}
\Phi_3(y),\quad \mbox{a.e.}\ \sigma(dy)du.
\end{equation}
Since $F$ is the kernel function of a PFSM, it satisfies
\begin{equation}\label{e:fourth-no-periodic3}
F(y,c(y)u) = b(y) F(y,u+a(y))+ d(y),\quad \mbox{a.e.}\ \sigma(dy)du,
\end{equation}
for some $c(y) >0\ (c(y)\neq 1)$, $b(y)\neq 0$, $a(y),d(y)\in\bbR$. Then, by replacing
$u$ by $c(y)u$ in (\ref{e:fourth-no-periodic2}) and by using
(\ref{e:fourth-no-periodic3}) and (\ref{e:fourth-no-periodic}), we get that
\begin{equation}\label{e:fourth-no-periodic4}
G^{C\setminus P}(\Phi_1(y),c(y)u) = B(y)G^{C\setminus P}(\Phi_1(y),u+A(y)) + D(y),\quad
\mbox{a.e.}\ \sigma(dy)du,
\end{equation}
for some $B(y)\neq 0$, $A(y),D(y)\in\bbR$. Since $\sigma(dy)$ is not a zero measure,
relation (\ref{e:fourth-no-periodic4}) contradicts the fact that $\Phi_1(y)\in C\setminus
C_P$ in view of the definitions of the set $C_P$.

The last two statements of the theorem follow from the proof of Theorem 5.3 in Pipiras
and Taqqu \cite{pipiras:taqqu:2002d}, Theorem 10.1 in Pipiras and Taqqu
\cite{pipiras:taqqu:2002s} and Propositions \ref{p:P-vs-C_P} and \ref{p:L=C_L}.
 \ \ $\Box$
\end{proof}


\section{Example of a process of the ``fourth'' kind}

We provide here examples of the ``fourth'' kind processes $X_\alpha^{C\setminus P}$ in
the decomposition (\ref{e:into-four-final}) which are related to $S\alpha S$
sub-Gaussian, more generally, sub-stable processes.

Let $\{W(t)\}_{t\in\bbR}$ be a stationary process which has c{\` a}dl{\` a}g (that is,
right continuous and with limits from the left) paths, satisfies $E|W(t)|^\alpha<\infty$,
\begin{equation}\label{e:for-fourth}
    E|W(t) - W(s)|^\alpha \leq C|t-s|^{2p},\quad s,t\in\bbR,
\end{equation}
for some $p>0$, $P(|W(t)| < c) < 1$ for all $c>0$ and is ergodic. Let also $\Omega =
\{w\}$ be the space of c{\` a}dl{\` a}g functions on $\bbR$ and $P(dw)$ be the
probability measure corresponding to the process $W$. Define a $S\alpha S$ stationary
process
$$
Y_\alpha^{(1)}(t) = \int_{\Omega} F(w,t) M_\alpha(dw),
$$
where $F(w,t) = w(t)$ and $M_\alpha(dw)$ has the control measure $P(dw)$. The process
$Y_\alpha^{(1)}$ is well-defined since $E|W(t)|^\alpha <\infty$. When the probability
measure $P$ corresponds to a Gaussian, more generally stable process, the process
$Y_\alpha^{(1)}$ is called sub-Gaussian, more generally sub-stable (see Samorodnitsky and
Taqqu \cite{samorodnitsky:taqqu:1994book}). The Lamperti transformation of the process
$Y_\alpha^{(1)}$ leads to a $S\alpha S$ self-similar process
$$
Y_\alpha^{(2)}(t) = \int_{\Omega} |t|^H F(w,\ln|t|) M_\alpha(dw).
$$
The process $Y_\alpha^{(2)}$ does not have stationary increments. We can transform it to
a process with stationary increments by the following procedure. Let
$$
Y_\alpha^{(3)}(t) = \int_{\Omega}\int_\bbR |t+u|^H F(w,\ln|t+u|) M_\alpha(dw,du),
$$
where $M_\alpha(dw,du)$ has the control measure $P(dw)du$. The process $Y_\alpha^{(3)}$
is self-similar and also stationary (in the sense of generalized processes). We can
transform it to a self-similar stationary increments process through the usual ``infrared
correction'' transformation $Y_\alpha^{(3)}(t) - Y_\alpha^{(3)}(0)$, that is,
\begin{equation}\label{e:example-fourth}
X_\alpha(t) = \int_\Omega\int_\bbR \Big(|t+u|^H F(w,\ln|t+u|) - |u|^H F(w,\ln|u|)
\Big)M_\alpha(dw,du).
\end{equation}

Observe that the process $X_\alpha$ is a self-similar mixed moving average by
construction. By Lemma \ref{l:fourth-well-defined} below, it is well-defined when
$H<\min\{p,1\}$. Moreover, the process $X_\alpha$ is generated by a conservative flow.
Indeed, by setting $G(w,u) = |u|^\kappa F(w,\ln|u|)$, we have $c^{-\kappa} G(w,cu) =
G(\psi_c(w),u)$, $c>0$, where
$$
\psi_z:w(z),z\in\bbR \mapsto w(z+\ln c),z\in\bbR,
$$
is a measurable flow on $\Omega$. Since the process $W(t),t\in\bbR$, is stationary, the
flow $\{\psi_c\}_{c>0}$ is measure preserving. It is conservative because the measure $P$
on $\Omega$ is finite and therefore there can be no wandering set of positive measure. By
Lemma \ref{l:fourth-P-empty} below, the PFSM set $C_P$ associated with the kernel in the
representation (\ref{e:example-fourth}) is empty a.e. Hence, in view of Theorem
\ref{t:into-four-final}, the process $X_\alpha$ is an example of the ``fourth'' kind
process $X_\alpha^{C\setminus P}$ in the decomposition (\ref{e:into-four-final}). We
state this result in the following proposition.

\begin{proposition}\label{p:fourth-kind}
The process $X_\alpha$ defined by (\ref{e:example-fourth}) under the assumptions stated
above, is an example of the process $X_\alpha^{C\setminus P}$ in the decomposition
(\ref{e:into-four-final}).
\end{proposition}

\medskip

The following auxiliary lemma shows that the process $X_\alpha$ in
(\ref{e:example-fourth}) is well-defined.

\begin{lemma}\label{l:fourth-well-defined}
The process $X_\alpha$ in (\ref{e:example-fourth}) is well-defined for
$H\in(0,\min\{p,1\})$ and $\alpha\in(0,2)$ under the assumption (\ref{e:for-fourth}).
\end{lemma}

\begin{proof}
The result follows since, by using (\ref{e:for-fourth}) and stationarity of $W$,
$$
\int_\Omega \int_\bbR \Big| |t+u|^\kappa F(w,\ln|t+u|) - |u|^\kappa
F(w,\ln|u|)\Big|^\alpha P(dw) du
$$
$$
= \int_\bbR E \Big| |t+u|^\kappa W(\ln|t+u|) - |u|^\kappa W(\ln|u|)\Big|^\alpha  du
$$
$$
\le 2^\alpha \int_\bbR  |t+u|^{\kappa\alpha} E\Big| W(\ln|t+u|) - W(\ln|u|)\Big|^\alpha
du + 2^\alpha \int_\bbR E|W(\ln|u|)|^\alpha \Big| |t+u|^\kappa - |u|^\kappa \Big|^\alpha
du
$$
$$
\leq 2^{\alpha} C \int_\bbR |t+u|^{\kappa\alpha} \Big|\ln|t+u| -\ln|u| \Big|^{p\alpha} du
+ 2^\alpha C \int_\bbR \Big| |t+u|^\kappa - |u|^\kappa \Big|^\alpha du < \infty,
$$
when $\kappa \alpha - p\alpha +1 = (H-1/\alpha)\alpha - p\alpha+1 = \alpha(H - p) <0$ and
$H<1$. \ \ $\Box$
\end{proof}

\medskip
The following lemma was used to show that the process $X_\alpha$ defined by
(\ref{e:example-fourth}) does not have a PFSM component.

\begin{lemma}\label{l:fourth-P-empty}
If $C_P$ is the PFSM set (\ref{e:pfsm-set}) associated with the representation
(\ref{e:example-fourth}) of the process $X_\alpha$, then $C_P = \emptyset$ a.e.\ $P(dw)$.
\end{lemma}

\begin{proof}
By the definition of the set $C_P$ in (\ref{e:pfsm-set}), we have
\begin{equation}\label{e:set-C_P-fourth}
C_P = \Big\{w\in \Omega: \exists c\neq 1,a,b\neq 0,d: |cu|^\kappa w(\ln|cu|) = b
|u+a|^\kappa w(\ln|u+a|) + d,\ \forall u\Big\},
\end{equation}
where the ``a.e.\ $du$'' condition in (\ref{e:pfsm-set}) was replaced by the ``$\forall
u$'' condition because the functions $w$ are c{\` a}dl{\` a}g. We may suppose without
loss of generality that $c>1$ in (\ref{e:set-C_P-fourth}). (If $c<1$, by making the
change of variables $u+a=c^{-1}v$ and dividing both sides of the relation in
(\ref{e:set-C_P-fourth}) by $b$, we obtain the relation analogous to
(\ref{e:set-C_P-fourth}) where $c$ is replaced by $c^{-1}$.) We shall consider the cases
$\kappa>0$ and $\kappa\leq 0$ separately.

\medskip
{\it The case $\kappa>0$:} We first examine the case when $b\neq 1$ in
(\ref{e:set-C_P-fourth}). By using (\ref{e:pfsm-set-exp2}) in Proposition
\ref{p:pfsm-set-alt}, we can express the equation in (\ref{e:set-C_P-fourth}) as
\begin{equation}\label{e:w=b*w}
|cu+g|^\kappa w(\ln|cu+g|)+f = b \Big(|u+g|^\kappa w(\ln|u+g|) + f\Big),
\end{equation}
for some $c>1,b\neq 0,f,g\in\bbR$. Setting
\begin{equation}\label{e:w-tildew}
    \widetilde w(v) = e^{-\kappa v} \Big(|e^v+g|^\kappa w(\ln|e^v+g|) + f\Big)
\end{equation}
and $\widetilde c = \ln c >0$, we have from (\ref{e:w=b*w}) with $u=e^v$ that
\begin{equation}\label{e:tildew=b*tildew}
    \widetilde w(v + \widetilde c) = \widetilde b\widetilde w(v),\quad v\in\bbR,
\end{equation}
where $\widetilde b = bc^\kappa$. Observe also that, by making the change of variables $v
= \ln (e^u-g)$ in (\ref{e:w-tildew}) for large $v$, we have
\begin{equation}\label{e:tildew-w}
    w(u) = e^{-\kappa u} \Big((e^u-g)^\kappa \widetilde w(\ln(e^u-g)) - f\Big)
\end{equation}
for large $u$.

If $|\widetilde b|\leq 1$  in (\ref{e:tildew=b*tildew}), then $|\widetilde w(v)|$ is
bounded for large $v$. Indeed, if $|\widetilde b|=1$, then $|\widetilde w(v)|$ is
periodic with period $\widetilde c$ and, being c{\` a}dl{\` a}g, it has to be bounded. If
$|\widetilde b|<1$, then $|\widetilde w(v)|\to 0$ as $v\to\infty$ because $|\widetilde
w(v + n\widetilde c)| = |\widetilde b|^n |\widetilde w(v)|$ and $|\widetilde b|^n\to 0$
as $n\to\infty$. By (\ref{e:tildew-w}), since $\kappa>0$, we obtain that $|w(u)|$ is
bounded for large $u$ as well. By Lemma \ref{l:zero-probabilities}, $(i)$, below, the
$P$-probability of such $w$ is zero.

Suppose now that $|\widetilde b|>1$ in (\ref{e:tildew=b*tildew}). We have either $(i)$
$\widetilde w(v)=0$ for $v\in [0,\widetilde c]$, or $(ii)$ $\inf\{|\widetilde w(v)|:v\in
A\}>0$ for $A\subset[0,\widetilde c]$  of positive Lebesgue measure. In the case $(i)$,
(\ref{e:tildew=b*tildew}) implies that $\widetilde w(v)=0$ for all $v$ and hence, by
(\ref{e:tildew-w}), $w(u)=-fe^{-\kappa u}$ for large $u$. By Lemma
\ref{l:zero-probabilities}, $(i)$, below, the $P$-probability of such $w$ is zero.
Consider now the case $(ii)$. Since $|\widetilde b|>1$, we get that
$$
\inf\{|\widetilde w(v)|: v\in A+n\widetilde c\} \to\infty,\quad \mbox{as}\ n\to\infty.
$$
Using (\ref{e:w-tildew}), since $\kappa>0$ (and hence $fe^{-\kappa u}\to 0$ as
$u\to\infty$), this yields that
$$
\inf\{|w(\ln(e^v+g))|: v\in A+n\widetilde c\} \to\infty,\quad \mbox{as}\ n\to\infty.
$$
By Lemma \ref{l:zero-probabilities}, $(iii)$, below, the $P$-probability of such $w$ is
zero.

If $b=1$ in (\ref{e:set-C_P-fourth}), by  using (\ref{e:pfsm-set-exp1}) in Proposition
\ref{p:pfsm-set-alt}, we get
\begin{equation}\label{e:w=b*w-b=1} |cu+g|^\kappa
w(\ln|cu+g|) = |u+g|^\kappa w(\ln|u+g|) + d,
\end{equation}
for some $c>1,g,d\in\bbR$. Setting
\begin{equation}\label{e:w-tildew-b=1}
    \widetilde w(v) = |e^v+g|^\kappa w(\ln|e^v+g|)
\end{equation}
and $\widetilde c = \ln c >0$, we deduce from (\ref{e:w=b*w-b=1}) with $u=e^v$ that
\begin{equation}\label{e:tildew=b*tildew-b=1}
    \widetilde w(v + \widetilde c) = \widetilde w(v) + d,\quad v\in\bbR.
\end{equation}
The function $\widetilde w$ is bounded on $[0,\widetilde c]$ since it is c{\` a}dl{\` a}g
and in view of (\ref{e:tildew=b*tildew-b=1}), we get
\begin{equation}\label{e:tildew-leq-Cv}
    |\widetilde w(v)|\leq C|v|,
\end{equation}
for large $v$ and some constant $C = C(w)>0$. Substituting (\ref{e:w-tildew-b=1}) into
(\ref{e:tildew-leq-Cv}), and since $\kappa >0$, we get that $w(v)\to 0$ as $v\to\infty$.
By Lemma \ref{l:zero-probabilities}, $(i)$, below, the $P$-probability of such $w$ is
zero. Combining this with the analogous conclusion  when $b\neq 1$ above, we deduce that
$C_P=\emptyset$ a.e.\ $P(dw)$ when $\kappa>0$.

\bigskip

{\it The case $\kappa\leq 0$:} By using (\ref{e:pfsm-set-exp1}) in Proposition
\ref{p:pfsm-set-alt}, we express the equation in (\ref{e:set-C_P-fourth}) as
\begin{equation}\label{e:w=b*w-0}
|cu+g|^\kappa w(\ln|cu+g|) = b |u+g|^\kappa w(\ln|u+g|) + d,
\end{equation}
for some $c>1,b\neq 0,d,g\in\bbR$. When $d=0$, we can use here the argument in the case
$\kappa>0$ because the assumption $\kappa>0$ was used above only to ensure that the term
$e^{-\kappa v}f$ in (\ref{e:w-tildew}) is negligible for large $v$. Suppose then $d\neq
0$. We can rewrite (\ref{e:w=b*w-0}) as
\begin{equation}\label{e:tildew=b*tildew-0}
    \widetilde w(v + \widetilde c) = b\widetilde w(v) + d,\quad v\in\bbR,
\end{equation}
where $\widetilde c = \ln c>0$ and
\begin{equation}\label{e:w-tildew-0}
    \widetilde w(v) = |e^v+g|^\kappa w(\ln|e^v+g|),\quad v\in\bbR.
\end{equation}
It follows from (\ref{e:tildew=b*tildew-0}) that
\begin{equation}\label{e:tildew-n=tildew-0}
    \widetilde w(v+n\widetilde c) = \left\{ \begin{array}{ll}
    b^n \widetilde w(v) + d\frac{b^n-1}{b-1}, & \mbox{if}\ b\neq 1,\\
    \widetilde w(v) + dn, & \mbox{if}\ b=1.
    \end{array}\right.
\end{equation}
Observe also that, by making the change of variable $v = \ln(e^u-g)$ in
(\ref{e:w-tildew-0}) for large $v$, we get
\begin{equation}\label{e:tildew-w-0}
    w(u) = e^{-\kappa u} \widetilde w(\ln(e^u-g)),
\end{equation}
for large $u$. We now consider separately the cases $|b|<1$, $|b|>1$, $b=1$ and $b=-1$.

$(a)$ Consider first the case $|b|<1$. By using (\ref{e:tildew-n=tildew-0}), we have
$$
\sup_{v\in[0,\widetilde c)}\Big| \widetilde w(v+n\widetilde c) + \frac{d}{b-1}\Big| =
|b|^n \sup_{v\in[0,\widetilde c)} \Big|\widetilde w(v) + \frac{d}{b-1}\Big| \to 0,\quad
\mbox{as} \ n\to\infty.
$$
Hence,
$$
\sup_{v\in[n\widetilde c,(n+1)\widetilde c)}\Big| \widetilde w(v) + \frac{d}{b-1}\Big|
\to 0,\quad \mbox{as} \ n\to\infty,
$$
and
\begin{equation}\label{e:tildew-|h|<1-0}
    \widetilde w(v)\to -\frac{d}{b-1},\quad \mbox{as}\ v\to\infty,
\end{equation}
or equivalently, by using (\ref{e:tildew-w-0}),
\begin{equation}\label{e:tildew-|h|<1-0,w}
    e^{\kappa u}w(u)\to -\frac{d}{b-1}\neq 0,\quad \mbox{as}\ u\to\infty.
\end{equation}
When $\kappa< 0$, relation (\ref{e:tildew-|h|<1-0,w}) implies that $|w(u)|\geq \epsilon
e^{-\kappa u}$ for large  $u$ and some constant $\epsilon>0$, that is, $|w(u)|$ is
unbounded for large $u$. When $\kappa=0$, we get that $|w(u)|$ is bounded for large $u$.
By Lemma \ref{l:zero-probabilities}, $(i)$--$(ii)$, below, the $P$-probability of such
$w$ in either case and hence those $w$ satisfying (\ref{e:w=b*w-0}) with $|b|<1$ is zero.

$(b)$ Consider now the case $|b|>1$. Relation (\ref{e:tildew-n=tildew-0}) can be
expressed as
\begin{equation}\label{e:tildew-n=tildew-0,1}
\widetilde w(v+n\widetilde c)+ \frac{d}{b-1} =
    b^n \Big(\widetilde w(v) + \frac{d}{b-1}\Big).
\end{equation}
We have either $(i)$ $\widetilde w(v) = -d/(b-1)$ for all $v$, or $(ii)$ there is a set
$A\subset [0,\widetilde c]$ of positive Lebesgue measure such that
$$
\inf\Big\{ \Big| \widetilde w(v) + \frac{d}{b-1} \Big|: v\in A \Big\} > 0.
$$
In the case $(i)$, by using (\ref{e:tildew-w-0}), we get that $w(u) = -d e^{-\kappa
u}/(b-1)$ for large $u$. The $P$-probability of such $w$ is zero by Lemma
\ref{l:zero-probabilities}, $(i)$--$(ii)$, below. In the case $(ii)$, by using
(\ref{e:tildew-n=tildew-0,1}) and since $|b|>1$, we have
$$
\inf\Big\{ \Big| \widetilde w(v) + \frac{d}{b-1} \Big|: v\in A+n\widetilde c \Big\} \to
\infty, \quad \mbox{as}\ n\to\infty
$$
or, in view of (\ref{e:w-tildew-0}) and $\kappa\leq 0$,
$$
\inf\{ | w(\ln(e^v+g))|: v\in A+n\widetilde c \} \to \infty, \quad \mbox{as}\ n\to\infty.
$$
The $P$-probability of such $w$ is zero by Lemma \ref{l:zero-probabilities}, $(iii)$,
below.

$(c)$ When $b=-1$ in (\ref{e:tildew=b*tildew-0}), we have $b^{2n}=1$ and $b^{2n}-1=0$,
and hence relation (\ref{e:tildew-n=tildew-0}) implies that
\begin{equation}\label{e:tildew-periodic-2c}
\widetilde w(v + n2\widetilde c) = \widetilde w(v),\quad v\in\bbR.
\end{equation}
Consider first the case $\kappa<0$. We have either $(i)$ $\widetilde w(v) = 0$ for all
$v$, or $(ii)$ there is a set $A\subset [0,2\widetilde c]$ of positive Lebesgue measure
such that
\begin{equation}\label{e:inf>0}
\inf\{ | \widetilde w(v) |: v\in A \}> 0.
\end{equation}
Arguing as in part $(i)$ above, the $P$-probability of $w$ satisfying this part $(i)$ is
zero. In the case $(ii)$, relations (\ref{e:tildew-periodic-2c}) and (\ref{e:inf>0})
imply that
$$
\inf\{ | \widetilde w(v) |: v\in A + n2\widetilde c \} = \inf\{ | \widetilde w(v) |: v\in
A \}> 0.
$$
By using (\ref{e:w-tildew-0}), and since $\kappa<0$,
$$
\inf\{ | w(\ln(e^v+g))|: v\in A+n\widetilde c \} \to \infty, \quad \mbox{as}\ n\to\infty.
$$
The $P$-probability of such $w$ is zero by Lemma \ref{l:zero-probabilities}, $(iii)$,
below. Turning to the case $\kappa=0$, relation (\ref{e:tildew-periodic-2c}) shows that
$|\widetilde w(v)|$ is periodic and hence bounded, since it is c{\` a}dl{\` a}g. By using
(\ref{e:tildew-w-0}), since $\kappa=0$, $|w(u)|$ is bounded for large $u$ as well. By
Lemma \ref{l:zero-probabilities}, $(i)$, the $P$-probability of such $w$ and hence of
those $w$ satisfying (\ref{e:w=b*w-0}) with $b=-1$ is zero.

$(d)$ When $b=1$ in (\ref{e:tildew=b*tildew-0}), relation (\ref{e:tildew-n=tildew-0})
becomes $\widetilde w(v + n\widetilde c) = \widetilde w(v) + dn$, $v\in\bbR$. Consider
the cases $(i)$ $\widetilde w(v) = 0$ for $v\in[0,\widetilde c]$, and $(ii)$ there is a
set $A\subset [0,\widetilde c]$ of positive Lebesgue measure such that $\inf\{ |
\widetilde w(v) |: v\in A \}> 0$. Arguing as above, the $P$-probability of $w$ satisfying
$(i)$ is zero. In the case $(ii)$, since $|d|n\to\infty$ as $n\to\infty$, we get that
$\inf\{ |\widetilde w(v)|: v\in A+n\widetilde c \} \to \infty$ as $n\to\infty$. By using
(\ref{e:w-tildew-0}), we get again that
$$
\inf\{ | w(\ln(e^v+g))|: v\in A+n\widetilde c \} \to \infty, \quad \mbox{as}\ n\to\infty.
$$
The $P$-probability of such $w$ is zero by Lemma \ref{l:zero-probabilities}, $(iii)$,
below. Combining the results for $|b|<1$, $|b|>1$, $b=-1$ and $b=1$, we conclude that the
$P$-probability of $w$ satisfying (\ref{e:w=b*w-0}) is zero. In other words, $C_P =
\emptyset$ a.e.\ $P(dw)$ when $\kappa\leq 0$ as well.
 \ \ $\Box$
\end{proof}

\bigskip
The next result was used in the proof of Lemma \ref{l:fourth-P-empty} above. Consider a
function $w:\bbR\mapsto \bbR$. We say that the function $|w(u)|$, $u\in\bbR$, is {\it
ultimately unbounded} if there is a set $A=A(w)\subset [0,C]$ of positive Lebesgue
measure with a fixed constant $C$ such that
$$
\inf\{|w(u)|:u\in A+nC\} \to \infty,\quad \mbox{as}\ n\to\infty.
$$
We say that $|w(u)|$, $u\in\bbR$, is bounded for large $u$ if there is $N=N(w)$ such that
$|w(u)|\leq N$ for large enough $u$. Denote
$$
A_1 = \{w:|w(u)|\ \mbox{{\rm is bounded for large}}\ u\},
$$
$$
A_2 = \{w:|w(u)|\ \mbox{{\rm is ultimately unbounded}}\},
$$
$$
A_3 = \{w:|w(\ln(e^u+g))|\ \mbox{{\rm is ultimately unbounded}}\},
$$
where $g=g(w)\in\bbR$.

\begin{lemma}\label{l:zero-probabilities}
Under the assumptions on the process $W$ (and hence on the corresponding probability $P$)
stated in the beginning of the section and with the sets $A_1,A_2,A_3$ defined above, we
have
$$
(i)\ P(A_1)=0,\quad  (ii)\ P(A_2)=0,\quad  \mbox{and}\quad (iii)\ P(A_3)=0.
$$
\end{lemma}

\begin{proof}
To show $(i)$, observe that $P(A_1)\leq \sum_{n=1}^\infty P(B_n)$, where $B_n = \{w:
|w(u)|< n\ \mbox{for large}\ u\}$. It is enough to show that $P(B_n)=0$ for $n\geq 1$.
When $w\in B_n$, we have
$$
\frac{1}{T} \int_0^T 1_{\{|w(u)|< n\}} du \to 1,
$$
as $T\to\infty$. But by ergodicity and the assumption $P(|w(0)|<c)<1$ for any $c>0$, we
have
$$
\frac{1}{T} \int_0^T 1_{\{|w(u)|< n\}} du \to P(|w(0)|<n) < 1\quad \mbox{a.e.}\ P(dw).
$$
This implies that $P(B_n)=0$.

We now show $(ii)$. Let $w\in A_2$, and $A$ and $C$ be the set and the constant appearing
in the definition of ultimate unboundedness of $|w|$. Observe that $\int_0^T
|w(u)|^\alpha du = \sum_{k=1}^K \int_{(k-1)C}^{kC} |w(u)|^\alpha du$ when $T = KC$, and
$\int_{(k-1)C}^{kC} |w(u)|^\alpha du \geq Leb(A) (\inf\{|w(u)|: u\in A+(k-1)C\})^\alpha
\to \infty$ as $k\to\infty$, by the ultimate unboundedness of $w$. Then, for $w\in A_2$,
we have
\begin{equation}\label{e:ergodicity-alpha-no}
\frac{1}{T} \int_0^T |w(u)|^\alpha du \to \infty.
\end{equation}
However, by ergodicity and the assumption $E|w(0)|^\alpha<\infty$, we have
\begin{equation}\label{e:ergodicity-alpha}
\frac{1}{T} \int_0^T |w(u)|^\alpha du \to E(|w(0)|^\alpha < \infty \quad \mbox{a.e.}\
P(dw).
\end{equation}
This implies that $P(A_2)=0$.

Consider now part $(iii)$. When $w\in A_3$ and $u=u_0$ is large enough, we have
$$
\frac{1}{T} \int_{u_0}^T |w(\ln(e^u+g))|^\alpha du \to \infty.
$$
Making the change of variables $\ln(e^u+g)=v$, we obtain that
$$
\frac{1}{T} \int_{\ln(e^{u_0}+g)}^{\ln(e^T+g)} |w(v)|^\alpha \frac{e^v}{e^v-g} dv \to
\infty.
$$
It is easy to see that this implies (\ref{e:ergodicity-alpha-no}) when $w\in A_3$. In
view of (\ref{e:ergodicity-alpha}), we get $P(A_3)=0$.
 \  \ $\Box$
\end{proof}


\small


\noindent Vladas Pipiras \hfill  Murad S.\ Taqqu

\noindent Department of Statistics and Operations Research    \hfill Department of
Mathematics and Statistics

\noindent  University of North Carolina at Chapel Hill \hfill Boston University

\noindent  CB\#3260, New West \hfill 111 Cummington St.

\noindent  Chapel Hill, NC 27599, USA \hfill Boston, MA 02215, USA

\noindent {\it pipiras@email.unc.edu} \hfill {\it murad@math.bu.edu}

\end{document}